\documentclass[10pt,leqno,dvipdfmx]{amsart}

\usepackage{amsmath,amsmath,amssymb,amsthm,braket,stmaryrd}

\usepackage[abbrev]{amsrefs}
\usepackage{graphicx}

\usepackage{mathtools}
\mathtoolsset{showonlyrefs=true}

\usepackage[top=30mm, bottom=30mm, left=25mm, right=25mm]{geometry}

\newtheorem{theo}{Theorem}[section]
\newtheorem{lemm}[theo]{Lemma}
\newtheorem{corr}[theo]{Corollary}

\numberwithin{equation}{section}

\theoremstyle{definition}

\newtheorem{rema}[theo]{Remark}
\newtheorem{assu}[theo]{Assumption}

\newcommand{\ba}{\mathbf{a}}

\newcommand{\bff}{\mathbf{f}}
\newcommand{\bg}{\mathbf{g}}
\newcommand{\bh}{\mathbf{h}}
\newcommand{\bk}{\mathbf{k}}

\newcommand{\bn}{\mathbf{n}}

\newcommand{\bu}{\mathbf{u}}
\newcommand{\bv}{\mathbf{v}}
\newcommand{\bw}{\mathbf{w}}

\newcommand{\bB}{\mathbf{B}}

\newcommand{\bD}{\mathbf{D}}

\newcommand{\bF}{\mathbf{F}}
\newcommand{\bG}{\mathbf{G}}
\newcommand{\bH}{\mathbf{H}}
\newcommand{\bI}{\mathbf{I}}

\newcommand{\bT}{\mathbf{T}}
\newcommand{\bV}{\mathbf{V}}

\newcommand{\DV}{\mathrm{Div}\,}
\newcommand{\dv}{\mathrm{div}\,}

\newcommand{\BR}{\mathbb{R}}
\newcommand{\BC}{\mathbb{C}}

\newcommand{\BN}{\mathbb{N}}

\newcommand{\BJ}{\mathbb{J}}

\newcommand{\BT}{\mathbb{T}}

\newcommand{\CA}{\mathcal{A}}

\newcommand{\CC}{\mathcal{C}}

\newcommand{\CJ}{\mathcal{J}}

\newcommand{\CM}{\mathcal{M}}

\newcommand{\CR}{\mathcal{R}}
\newcommand{\CS}{\mathcal{S}}

\newcommand{\CU}{\mathcal{U}}

\newcommand{\fp}{\mathfrak{p}}

\newcommand{\pd}{\partial}

\newcommand{\lbrac}{\llbracket}
\newcommand{\rbrac}{\rrbracket}
\newcommand{\ov}{\overline}
\newcommand{\dx}{\,\mathrm{d}x}
\newcommand{\ds}{\,\mathrm{d}s}
\newcommand{\dt}{\,\mathrm{d}t}

\newcommand{\wh}{\widehat}
\newcommand{\wt}{\widetilde}

\begin{document}

\title[]{Global solvability of compressible--incompressible two-phase flows with phase transitions in bounded domains}

\author[]{Keiichi Watanabe}

\address{Global Center for Science and Engineering, Waseda University, 3-4-1 Ookubo, Shinjuku-ku, Tokyo, 169-8555, Japan}		


\email{keiichi-watanabe@akane.waseda.jp}

\thanks{This research was partly supported by JSPS Grant-in-aid for Research Activity Start-up Grant Number 20K22311 and Waseda University Grant for Special Research Projects.}


\date{}

\maketitle

\begin{abstract}
Consider a free boundary problem of compressible-incompressible two-phase flows with surface tension and phase transition in bounded domains $\Omega_{t +}, \Omega_{t -} \subset \BR^N$, $N \ge 2$, where the domains are separated by a sharp compact interface $\Gamma_t \subset \BR^{N - 1}$. We prove a global in time unique existence theorem for such free boundary problem under the assumption that the initial data are sufficiently small and the initial domain of the incompressible fluid is close to a ball. In particular, we obtain the solution in the maximal $L_p - L_q$-regularity class with $2 < p <\infty$ and $N < q < \infty$ and exponential stability of the corresponding analytic semigroup on the infinite time interval.
\end{abstract}

\section{Introduction}
\noindent
The purpose of this paper is to prove the global solvability of the free boundary problem of
compressible--incompressible two-phase flows with phase transitions in bounded domains. Two fluids are separated by a free boundary $\Gamma_t$ and a surface tension and phase transitions are taken into account. Our problem is formulated as follows: Let $\Omega$ be a bounded domain in $N$-dimensional Euclidean space $\BR^N$, $N \ge 2$, surrounded by a smooth boundary $\Gamma_+$. For $t \ge 0$, the hypersurface $\Gamma_t$ represents a sharp moving interface that separates $\Omega$ into $\Omega_{t +}$ and $\Omega_{t -}$ such that $\Omega \setminus \Gamma_t = \Omega_{t +} \cup \Omega_{t -}$. Through this article, we suppose that the two fluids are immiscible and no boundary contact occurs. Let $\dot \Omega_t = \Omega_{t +} \cup \Omega_{t -}$ and for any function $f$ defined on $\dot \Omega_t$, we write $f_\pm = f \vert_{\Omega_{t \pm}}$. We consider the following Cauchy problem:
\begin{align}
\label{eq-1.1}
&\left\{\begin{aligned}
\pd_t\varrho_+ + \dv (\varrho_+ \bv_+) & = 0 & \qquad & \text{in $\Omega_{t +}, \enskip t > 0$},\\
\varrho_+ (\pd_t \bv_+ + (\bv_+ \cdot \nabla) \bv_+) - \DV \BT_+ & = 0 & \qquad & \text{in $\Omega_{t +}, \enskip t > 0$},\\
\dv \bv_- & = 0 & \qquad & \text{in $\Omega_{t -}, \enskip t > 0$}, \\
\varrho_- (\pd_t \bv_- + (\bv_- \cdot \nabla) \bv_-) - \DV \BT_- & = 0 & \qquad & \text{in $\Omega_{t -}, \enskip t > 0$}	
\end{aligned}\right.
\intertext{with the interfacial boundary conditions on $\Gamma_t$ $(t > 0)$:}
\label{eq-1.2}
&\left\{\begin{aligned}
& V_{\Gamma_t} = \bv_{\Gamma_t} \cdot \bn_t = \cfrac{\lbrac \varrho \bv \rbrac \cdot \bn_t}{\lbrac \varrho \rbrac},\\
& \lbrac \bv \rbrac = j \bigg \lbrac \frac{1}{\varrho} \bigg \rbrac \bn_t, 
\quad j \lbrac \bv \rbrac - \lbrac \BT \rbrac \bn_t = - \sigma H_{\Gamma_t} \bn_t,\\
& \lbrac \psi \rbrac + \cfrac{j^2}{2} \bigg \lbrac \cfrac{1}{\varrho^2} \bigg \rbrac 
- \bigg \lbrac \cfrac{1}{\varrho} (\BT \bn_t \cdot \bn_t) \bigg \rbrac = 0,\\
& (\nabla \varrho_+) \cdot \bn_t \rvert_+ = 0,	
\end{aligned}\right.	
\end{align}
and the homogeneous Dirichlet boundary conditions on $\Gamma_+$:
\begin{equation}
\label{cond-1.3}
\begin{aligned}
\bv_+ = 0 , \quad \nabla \varrho_+ \cdot \bn_+ = 0 \quad \text{on $\Gamma_+$},	
\end{aligned}	
\end{equation}
and the initial conditions:
\begin{equation}
\label{cond-1.4}
\begin{aligned}
(\varrho_+, \bv_+)\rvert_{t = 0} = (\rho_{* +} + \rho_{0 +}, \bv_{0 +}) \quad \text{in $\Omega_{0 +}$}, \quad
\bv_-\rvert_{t = 0} = \bv_{0 -} \quad \text{in $\Omega_{0 -}$},	\quad 
\Gamma_t \vert_{t = 0} = \Gamma_0,
\end{aligned}	
\end{equation}
where $\varrho_+$ and $\varrho_- := \rho_{* -}$ are the densities, $\bv_\pm$ the velocity fields, $\psi_\pm$ the Helmholtz free energy functions, and $\rho_{* \pm}$ are positive constants, $H_{\Gamma_t}$ the $(N - 1)$-times mean curvature of $\Gamma_t$, $\sigma$ a positive constant describing the coefficient of the surface tension, $V_t$ the velocity of evolution of $\Gamma_t$ with respect to $\bn_t$, $\bv_{\Gamma_t}$ the interfacial velocity, $\bn_t$ the outer unit normal to $\Gamma_t$ pointed from $\Omega_{t +}$ to $\Omega_{t -}$, and $\bn_+$ the outer unit normal to $\Gamma_+$. Here, $j = \varrho_+ (\bv_+ - \bv_\Gamma) \cdot \bn_t = \varrho_- (\bv_- - \bv_\Gamma)\cdot\bn_t$ is the phase flux and $\BT_\pm$ are the Stress tensors defined by
\begin{equation*}
\begin{split}
\BT_+ = & \mu_+ \bD(\bv_+) + (\nu_+ - \mu_+) \dv \bv_+ \bI - \fp_+ \bI + \bigg( \frac{\kappa_+}{2} \lvert \nabla \varrho_+ \rvert^2
+ \kappa_+ \varrho_+ \Delta \varrho_+ \bigg) \bI - \kappa_+ \nabla \varrho_+ \otimes \nabla \varrho_+,\\
\BT_- = & \mu_- \bD(\bv_-) - \fp_- \bI,	
\end{split}	
\end{equation*}
where $\fp_\pm$ are the pressure fields. Here, for any vector fields $\ba = (a_1, \dots , a_N)$, the deformation tensor $\bD (\ba)$ is defined by $\bD(\ba) = \nabla \ba + {}^\top\! (\nabla \ba)$ whose $(j, k)$th components are given by $\pd a_k + \pd_k a_j$. Besides,
\begin{equation*}
\bigg( \frac{\kappa_+}{2} \lvert \nabla \varrho_+ \rvert^2
+ \kappa_+ \varrho_+ \Delta \varrho_+ \bigg) \bI - \kappa_+ \nabla \varrho_+ \otimes \nabla \varrho_+
\end{equation*}
is called the Korteweg tensor, which essentially ensures an additional smoothing for the density. In this article, the coefficients $\mu_+$, $\nu_+$, $\kappa_+$, and $\mu_-$ are assumed to be constants satisfying
\begin{equation*}
\mu_+ > 0, \quad \nu_+ \ge \frac{N - 2}{N} \mu_+, \quad \kappa_+ > 0, \quad \mu_- > 0.
\end{equation*}
Notice that if phase transitions occur on the moving interface $\Gamma_t$, the phase flux $j$ should be taken arbitrary. Furthermore, the jump of a quantity $\bg(x, t)$ defined on $\dot \Omega_t$ across the interface $\Gamma_t$ is defined by 
\begin{equation*}
\lbrac \bg \rbrac (x_0) := \lim_{\delta \to 0 +} \left(\bg (x_0 + \delta \bn_t(x_0) ) - \bg (x_0 - \delta \bn_t (x_0)) \right)	
\end{equation*}
for all $x_0 \in \Gamma_t$, where $\bn_t (x_0)$ is the outer unit normal to $\Gamma_t$ at $x_0$. In addition, we adopt the notations $\bg \rvert_\pm (x_0) = \lim_{\delta \to 0 +} \bg (x_0 \mp \delta \bn_t(x_0))$ for all $x_0 \in \Gamma_t$. \par
Pr{\"u}ss et al. \cite{PS12,PSW14,PSSS12,PS16} and Shimizu and Yagi~\cite{SY15,SY17} studied the thermodynamically consistent model of the incompressible and incompressible two-phase flows with phase transitions. In particular, Pr{\"u}ss, Shimizu, and Wilke~\cite{PSW14} proved the stability of the equilibria of the problem. On the other hand, the compressible and incompressible two-phase flows with phase transitions was studied by Shibata~\cite{S16}. However, his result includes the derivative loss in the nonlinear term with respect to $\varrho_+$ due to the kinetic equation: $\bv_{\Gamma_t} \cdot \bn_t = \lbrac \varrho \bv \rbrac \cdot \bn_t / \lbrac \varrho \rbrac$. Namely, we can not prove the local solvability of the problem based on his result. To overcome this difficulty, the new model using the Navier--Stokes--Korteweg equations was proposed by the author~\cite{W18}. The new model is an extension of the Navier--Stokes--Fourier equations, and the physical consistency was discussed in~\cite[Sec.~2]{W18}. Here, the second and third conditions of \eqref{eq-1.2} stand for the conservation laws of mass and momentum on $\Gamma_t$, respectively. Besides, the condition $(\nabla \varrho_+) \cdot \bn_t \vert_+ = 0$ not only guarantees the generalized Gibbs-Thomson law
\begin{equation*}
\lbrac \psi \rbrac + \cfrac{j^2}{2} \bigg \lbrac \cfrac{1}{\varrho^2} \bigg \rbrac 
- \bigg \lbrac \cfrac{1}{\varrho} (\BT \bn_t \cdot \bn_t) \bigg \rbrac = 0 \qquad \text{on $\Gamma_t$}
\end{equation*}
but also implies the interstitial working: $(\kappa_+ \varrho_+ \dv \bv_+) \nabla \varrho_+$ vanishes in the normal direction of the interface $\Gamma_t$.
Recently, the local solvability of the problem was showed in the previous paper~\cite{W20}. For further historical review or physical backgrounds of our model, the readers may consult the introduction in \cite{PSSS12,PS16,W18,W20} and references therein. \par 
Now, let us formulate the problem. Let $B_R = \{y \in \BR^N \,\colon\, \lvert y \rvert < R\}$ and $S_R = \{y \in \BR^N \,\colon\, \lvert y \rvert = R\}$. In this paper, we suppose the following.
\begin{assu}
\label{asp-1}
Let $\lvert D \rvert$ denote the Lebesgue measure of a Lebesgue measurable set $D$ in $\BR^N$. We assume the following conditions:
\begin{enumerate}
\item It holds $\lvert \Omega_{0 -} \rvert = \lvert B_R \rvert = R^N \omega_N \slash N$, where $\omega_N$ is the area of $S_1$.
\item The barycenter point of $\Omega_{0 -}$ is suited at the origin, i.e.,
\begin{equation*}
\frac{1}{\lvert \Omega_{0 -} \rvert} \int_{\Omega_{0 -}} \rho_{* -} x \dx = 0.
\end{equation*}
\item The initial free surface $\Gamma_0$ is a normal perturbation of $S_R$ given by
\begin{equation*}
\Gamma_0 = \{x = y + h_0 (y) (y \slash \lvert y \rvert) \,\colon\, y \in S_R\} = \{x = (1 + R^{- 1} h_0 (y)) y \,\colon\, y \in S_R\}
\end{equation*}
with given small function $h_0$ defined on $S_R$.
\end{enumerate}
\end{assu}
Let the free boundary $\Gamma_t$ be given by
\begin{equation*}
\begin{split}
\Gamma_t & = \{x = y + h (y, t) (y \slash \lvert y \rvert) + \xi (t) \,\colon\, y \in S_R\} \\
& = \{x = (1 + R^{- 1} h(y, t)) y + \xi (t) \,\colon\, y \in S_R\},
\end{split}
\end{equation*}
where $h (y, t)$ is an unknown function such that $h (y, 0) = h_0 (y)$ for $y \in S_R$. Here, $\xi (t)$ denotes the barycenter point of the \textit{incompressible} domain $\Omega_{t -}$ defined by
\begin{equation*}
\xi (t) = \frac{1}{\lvert \Omega_{t -} \rvert} \int_{\Omega_{t -}} \rho_{* -} x \dx,
\end{equation*}
which is also an unknown function. Here, we have $\xi (0) = 0$ due to Assumption~\ref{asp-1} (2). Since it holds $\dv \bv_- = 0$ in $\Omega_{t -}$, the transport theorem implies
\begin{equation*}
\frac{\mathrm{d} \lvert \Omega_{t -} \rvert}{\dt} = \int_{\Omega_{t -}} \dv \bv \dx = 0,
\end{equation*}
which means that $\lvert \Omega_{t -} \rvert$ is independent of $t$, i.e., it holds $\lvert \Omega_{t -} \rvert = \lvert \Omega_{0 -} \rvert$. Hence, from the transport theorem, we have
\begin{equation*}
\label{xi-derivative}
\frac{\mathrm{d}}{\dt} \xi (t) = \frac{\rho_{* -}}{\lvert \Omega_{0 -} \rvert} \int_{\Omega_{t -}} \bv (x, t) \dx = \frac{\rho_{* -}}{\lvert B_R \rvert} \int_{\Omega_{t -}} \bv (x, t) \dx.
\end{equation*}
Furthermore, the transport theorem also yields an important formula
\begin{equation}
\label{transprot-omegaplus}
\frac{\mathrm{d} \lvert \Omega_{t +} \rvert}{\dt} = \int_{\Omega_{t +}} \dv \bv_+ \dx = 0
\end{equation} 
because $\lvert \Omega \rvert$ and $\lvert \Omega_{t -} \rvert$ are independent of $t$ so that $\lvert \Omega_{t + } \rvert$ is independent of $t$ as well. For given function $h(y, t)$, let $H_h (y, t)$ be a solution to the Dirichlet problem: $(1 - \Delta) H_h = 0$ in $\BR^N \setminus S_R$, $H_h = h$ on $S_R$. From the $K$-method in real interpolation theory \cite{L18}, we have
\begin{equation*}
\begin{aligned}
C_1 \lVert H_h (\cdot, t) \rVert_{H^k_q (\BR^N)} & \le \lVert h (\cdot, t) \rVert_{W^{k - 1 \slash q}_q (S_R)} \le C_2 \lVert H_h (\cdot, t) \rVert_{H^k_q (\BR^N)} & \quad & \text{for $k = 1,2,3$}, \\
\lVert \pd_t H_h (\cdot, t) \rVert_{H^k_q (\BR^N)} & \le \lVert \pd_t h (\cdot, t) \rVert_{W^{k - 1 \slash q}_q (S_R)} \le C_2 \lVert \pd_t H_h (\cdot, t) \rVert_{H^k_q (\BR^N)} & \quad & \text{for $k = 1,2$}.
\end{aligned}
\end{equation*}
We may assume that there exists a small number $r > 0$ such that $B_{R + 3 r} \subset \Omega$. Let $\varphi \in C^\infty (\BR^N)$ be a cutoff function that equals one for $x \in B_{R + r}$ and zero for $x \notin B_{R + 2 r}$. Let $\Phi (y, t) = y + \varphi (y) (R^{- 1} H_h (y, t) y + \xi (t))$. Notice that $\Phi (y, t) = y + R^{- 1} H_h (y, t) y + \xi (t)$ for $y \in B_R$. Setting $\Psi (y, t) = \varphi (y) (R^{- 1} H_h (y, t) + \xi (t))$, we assume that
\begin{equation}
\label{cond-psi}
\sup_{t \in (0, T)} \lVert \Psi (\cdot, t) \rVert_{H^1_\infty (\BR^N)} \le \delta
\end{equation}
with some small constant $\delta > 0$. In the following, we choose $\delta > 0$ so small that the map $y \mapsto x = \Phi (y, t)$ is bijective from $\Omega$ onto itself. In fact, for any $y_1$ and $y_2$, it holds
\begin{equation*}
\lvert \Phi (y_1, t) - \Phi (y_2, t) \rvert \ge \lvert y_1 - y_2 \rvert - \sup_{t \in (0, T)} \lVert \nabla \Phi (\cdot, t) \rVert_{H^1_\infty (\BR^N)} \lvert y_1 - y_2 \rvert \ge (1 - \delta) \lvert y_1 - y_2 \rvert,
\end{equation*}
which implies the injectivety of the map $x = \Phi (y, t)$ for any $t \ge 0$ provided that $\delta \in (0, 1)$. Furthermore, by the inverse mapping theorem, the map $x = \Phi (y, t)$ is surgective from $\Omega$ onto itself since $x = \Phi (y, t) = y$ for $y \in \Omega \setminus B_{R + 2 r}$. Let
\begin{equation*}
\begin{split}
\Omega_{t +} & = \{x = \Phi (y, t) = y + \varphi (y) (R^{- 1} H_h (y, t) y + \xi (t)) \,\colon\, y \in \Omega \setminus \ov{B_R}\}, \\
\Omega_{t -} & = \{x = \Phi (y, t) = y + R^{- 1} H_h (y, t) y + \xi (t) \,\colon\, y \in B_R\}, \\
\Gamma_t & = \{x = y + R^{- 1} h (y, t) y + \xi (t) \,\colon\, y \in S_R\}.
\end{split}
\end{equation*}
Here, $R^{- 1} y$ is the unit outer normal to $S_R$ for $y \in S_R$. \par
If one deals with the global existence issue of the free boundary problem of the Navier--Stokes equations with surface tension in a bounded domain, it is known that spectral analysis of the Stokes operator and the Laplace-Beltrami operator are crucial. To this end, we follow the approach due to Shibata~\cite{S18} in order to fixed the free boundary, where the corresponding transformation is given by $x = \Phi (y, t)$. The essential point of his approach is that an eigenvalue of the principal linearization does not appear on $\BC_+ = \{\lambda \in \BC \,\colon\, \mathrm{Re}\, \lambda \ge 0\}$, which yields the exponential stability of solutions as follows from the standard semigroup theory. Here, the similar approach was also used for the incompressible-incompressible two-phase flow case~\cite{EKS,SS20}. In our case, however, the domain of $\Omega_{t +}$ is occupied by the compressible fluid so that further dedication is required, where the relation \eqref{transprot-omegaplus} becomes crucial. The details will be explained in Section~\ref{sec-derivation}. \par
Let $\Phi_0 (y) = y + \varphi (y) (R^{- 1} H_{h_0} (y) y)$, where $H_{h_0}$ is a unique solution of the Dirichelt problem: $(1 - \Delta) H_{h_0} = 0$ in $\BR^N \setminus S_R$ and $H_{h_0} = h_0$ on $S_R$. In the following, we set $\rho_{0 +} = \varrho_{0 +} \circ \Phi_0$ and $\bu_{0 \pm} = \bv_{0 \pm} \circ \Phi_0$,
where $\varrho_{0 +}$ and $\bv_{0 \pm}$ are initial data \eqref{cond-1.4}. For functions $\varrho_+$, $\bv_\pm$, and $\fp_-$ satisfying the system \eqref{eq-1.1}--\eqref{cond-1.3}, we set
\begin{align*}
\rho_+ (y, t) = \varrho_+ \circ \Phi - \rho_{* +}, \quad \bu_\pm (y, t) = \bv_\pm \circ \Phi, \quad \pi_- (y, t) = \fp_- \circ \Phi - \frac{\sigma (N - 1)}{R}.
\end{align*}
Then the fixed boundary system associated with~\eqref{eq-1.1},~\eqref{eq-1.2}, \eqref{cond-1.3}, and~\eqref{cond-1.4} can be read as the following:
\begin{align}
\label{eq-1.5}
\left\{\begin{aligned}
\pd_t \rho_+ + \rho_{* +} \dv \bu_+ & = f_M (\rho_+, \bu_+, h) &\enskip &\text{ in $\Omega_+ \times(0, T)$}, \\
\rho_{* -} \dv \bu_- = f_d (\bu_-, h) & = \rho_{* -} \dv \bff_d (\bu_-, h) &\enskip &\text{ in $\Omega_- \times (0, T)$}, \\
\rho_{* +} \pd_t \bu_+ - \DV \bT_+ (\rho_+, \bu_+) & = \bff_+ (\rho_+, \bu_+, h) &\enskip &\text{ in $\Omega_+ \times (0, T)$}, \\
\rho_{* -} \pd_t \bu_- - \DV \bT_- (\bu_-, \pi_-) & = \bff_- (\bu_-, h) &\enskip &\text{ in $\Omega_- \times (0, T)$}, \\
\pd_t h - \frac{1}{\rho_{* -} - \rho_{* +}} \lbrac \langle \rho_* \bu, \bn \rangle \rbrac + \CM \bu & = d (\rho_+, \bu_+, \bu_-, h) &\enskip &\text{ on $S_R \times (0, T)$}, \\
\bB (\rho_+, \bu_+, \bu_-, \pi_-, h) & = \bG (\rho_+, \bu_+, \bu_-, h), &\enskip &\text{ on $S_R \times (0, T)$}, \\
\bu_+ = 0, \qquad \langle \nabla \rho_+, \bn_+ \rangle & = 0 &\enskip &\text{ on $\Gamma_+ \times (0, T)$}, \\
(\rho_+, \bu, h) \rvert_{t = 0} & = (\rho_{0 +}, \bu_0, h_0) &\enskip &\text{ on $\Omega_+ \times \dot \Omega \times S_R$},
\end{aligned}\right.	
\end{align}
where $\bn$ denotes the the outer unit normal to $S_R$ pointed from $\Omega_+$ into $\Omega_-$ and we have set
\begin{equation*}
\begin{split}
& \CM \bu = \frac{\rho_{* -}}{\lvert B_R \rvert} \bn \cdot \int_{\Omega_-} \bu_- \,\mathrm{d}y - \frac{\rho_{* +} \varphi_0}{\rho_{* -} - \rho_{* +}}\int_{\Omega_+} \dv \bu_+ \,\mathrm{d} y, \\
& \text{and} \qquad \frac{1}{\rho_{* -} - \rho_{* +}} \lbrac \langle \rho_* \bu, \bn \rangle \rbrac = \frac{1}{\rho_{* -} - \rho_{* +}} \Big(\langle \rho_{* -} \bu_-, \bn \rangle \rvert_- - \langle \rho_{* +} \bu_+, \bn \rangle \rvert_+\Big)
\end{split}
\end{equation*}
with $\omega = y \slash \lvert y \rvert \in S_1$ and $\varphi_0 = \lvert S_R \rvert^{- 1 \slash 2}$. By abuse of notation, here and in the following, we may write
\begin{equation*}
\bu_0 = \begin{cases}
\bu_{0 +} & \text{in $\Omega_+$}, \\
\bu_{0 -} & \text{in $\Omega_-$},	
\end{cases}
\qquad
\bu = \begin{cases}
\bu_+ & \text{in $\Omega_+$}, \\
\bu_- & \text{in $\Omega_-$},
\end{cases}
\qquad
\rho_* \bu = \begin{cases}
\rho_{* +} \bu_+ & \text{in $\Omega_+$}, \\
\rho_{* -} \bu_- & \text{in $\Omega_-$},
\end{cases}	
\end{equation*}
and we let $\bT_+$ and $\bT_-$ be ``linearized'' stress tensors defined by
\begin{align*}
\bT_+ (\bu_+, \rho_+) & := \mu_+ \bD (\bu_+) + (\nu - \mu) (\dv \bu_+) \bI + (- \gamma_{* +} + \rho_{* +} \kappa_+ \Delta) \rho_+ \bI, \\
\bT_- (\bu_-, \pi_-) & := \mu_- \bD(\bu_-) - \pi_- \bI
\end{align*}
and $\bB (\rho_+, \bu_+, \bu_-, \pi_-, h) = \bG (\rho_+, \bu_+, \bu_-, h)$ stands for the following interface conditions on $S_R \times (0, T)$:
\begin{align}
\label{def-B}
\left\{\begin{aligned}
\Pi_\bn (\mu_- \bD (\bu_-) \bn) \rvert_- - \Pi_\bn (\mu_+ \bD(\bu_+) \bn) \rvert_+ & = g (\rho_+, \bu_+, \bu_-, h), \\
\langle \bT_- (\bu_-, \pi_-) \bn, \bn \rangle \rvert_-
- \langle \bT_+ (\bu_+, \rho_+) \bn, \bn \rangle \rvert_+ + \sigma \CA_{S_R} h + \sum_{j = 1}^N (h, \varphi_j)_{S_R} \varphi_j & = f^+_B (\rho_+, \bu_+, \bu_-, h), \\
\frac{1}{\rho_{* -}} \langle \bT_- (\bu_-, \pi_-) \bn, \bn \rangle \rvert_- - \frac{1}{\rho_{* +}} \langle \bT_+ (\bu_+, \rho_+) \bn, \bn \rangle \rvert_+ + \frac{1}{\rho_{* -}} \sum_{j = 1}^N (h, \varphi_j)_{S_R} \varphi_j & = f^-_B (\rho_+, \bu_+, \bu_-, h),\\
\Pi_\bn \bu_- \rvert_- - \Pi_\bn \bu_+ \rvert_+ & = \bh(\bu_+, \bu_-, h), \\
\langle \nabla \rho_+, \bn \rangle \rvert_+ & = k_-(\rho_+, h),		
\end{aligned}\right.
\end{align}
where $\gamma_{* +} = \fp_+' (\rho_{* +})$, $\CA_{S_R} = - \Delta_{S_R} - (N - 1) \slash R^2$, and $\Pi_\bn \ba = \ba - \langle \ba, \bn \rangle \bn$ for any $N$ vector $\ba$. Furthermore, $\varphi_\ell$ ($\ell = 1, \dots, N$) denote the spherical harmonics of degree 1 on $S_R$, where $(\varphi_k,\varphi_\ell)_{S_R} = \delta_{k \ell}$. The right-hand members of~\eqref{eq-1.5} and \eqref{def-B} stands for the nonlinear terms that are independent of $\pi_-$, which will be explained in the next section. \par
Before stating our main results, we finally introduce some technical assumptions. 
\begin{assu}
\label{asp-2}
It holds $\rho_{* +} \ne \rho_{* -}$. The coefficients $\rho_{* +}$, $\mu_+$, $\nu_+$, $\kappa_+$, and $\mu_-$ satisfy
\begin{equation}
\label{cond-coeffi}
\bigg(\frac{\mu_+ + \nu_+}{2 \rho_{* +}^2 \kappa_+}\bigg)^2 \ne \frac{1}{\rho_{* +} \kappa_+}, \qquad \rho_{* +}^3 \kappa_+ = \mu_+ \nu_+.
\end{equation}
We further assume the following properties:
\begin{enumerate}
\item The pressure field $\fp_+ (\rho_+)$ is a $C^2$-function defined on $\rho_{* +} \slash 3 \le \rho_+ \le 3 \rho_{* +}$ such that $0 < \fp_+' (\rho_+) \le \pi^*$ with some positive constant $\pi^*$ for any $\rho_{* +} \slash 3 \le \rho_+ \le 3 \rho_{* +}$.
\item The Helmholtz free energy $\psi_+ (\rho_+, \lvert \rho_+ \rvert^2)$ is a $C^2$-function defined on $(\rho_{* +} \slash 3, 3 \rho_{* +}) \times [0, \infty)$ such that $0 \le \pd_{\varrho_+} \psi_+ (\rho_+, \lvert \rho_+ \rvert^2) \le \psi^*$ with some positive constant $\psi^*$ for any $\rho_{* +} \slash 3 \le \rho_+ \le 3 \rho_{* +}$. Besides, we assume that $\pd_{\varrho_+} \psi_+ (\rho_{* +}, 0) = 0$.
\item There exists positive constants $\pi_{* \pm}$ such that
\begin{equation*}
\psi_- (\rho_{* -}) - \psi_+ (\rho_{* +}, 0) = \frac{\pi_{* +}}{\rho_{* +}} - \frac{\pi_{* -}}{\rho_{* -}}, \qquad \pi_{* -} - \pi_{* +} = \sigma H_{S_R},
\end{equation*}
which stands for the Gibbs-Thomson condition and the Young-Laplace law, respectively. Especially, $\pi_{* +}$ is given by $\pi_{* +} = \fp_+ (\rho_{* +})$.
\end{enumerate}
\end{assu}
\begin{rema}
The conditions \eqref{cond-coeffi} are imposed to avoid multiple roots of the characteristic equation arising in the model problems in the half space and the whole space with flat interface. In fact, in those cases, applying the partial Fourier transform to the generalized resolvent problem yields the ODEs with respect to $x_N$, and the solution formula is obtained by the inverse partial Fourier transform, see \cite[Sect.~2]{Sai20} and \cite[Sect.~4]{W18}. The condition \eqref{cond-coeffi} expect to be removed by employing the similar argument due to Saito \cite{Sai20}.
\end{rema}
We now state our main result of this article. To this end, we record the necessary compatibility conditions for the given function $h_0$. According to Assumption~\ref{asp-1}, it follows that $h_0$ should satisfy the following conditions:
\begin{equation*}
\begin{split}
\frac{R^N \omega_N}{N} & = \int_{\Omega_{0 -}} \dx = \int_{\lvert \omega \rvert = 1} \int_0^{R + h_0 (R \omega)} r^{N - 1} \,\mathrm{d} r \, \mathrm{d} \omega = \int_{\lvert \omega \rvert = 1} \frac{1}{N} (R + h_0 (R \omega))^N \,\mathrm{d} \omega, \\
0 & = \int_{\Omega_{0 -}} x \dx = \int_{\lvert \omega \rvert = 1} \int_0^{R + h_0 (R \omega)} r^N \omega \,\mathrm{d} r \, \mathrm{d} \omega = \int_{\lvert \omega \rvert = 1} \frac{1}{N + 1} (R + h_0 (R \omega))^{N + 1} \omega \,\mathrm{d} \omega,
\end{split}
\end{equation*}
where $\mathrm{d} \omega$ denotes the surface element of $S_R$. Namely, we have the compatibility condition for $h_0$ as follows:
\begin{equation}
\label{compati-h0}
\begin{split}
\sum_{k = 1}^N {}_N \mathsf{C}_k \int_{S_R} (R^{- 1} h_0 (y)) \,\mathrm{d} \omega & = 0, \\ \sum_{k = 1}^{N + 1} {}_{N + 1} \mathsf{C}_k \int_{S_R} (R^{- 1} h_0 (y))^k y \,\mathrm{d} \omega & = 0,
\end{split}
\end{equation}
where ${}_j \mathsf{C}_k = \frac{j !}{k ! (j - k) !}$, $j \in \{N, N + 1\}$, are the binomial coefficients. Finally, we set
\begin{equation*}
(\rho_+, \bu_+, \bu_-, \pi_-, h) \in \CS_{p, q} (0, T) \quad \Leftrightarrow \quad \left\{\begin{split}
\rho_+ & \in H^1_p (0, T; H^1_q (\Omega_+)) \cap L_p (0, T; H^3_q (\Omega_+)), \\	
\bu_+ & \in H^1_p (0, T; L_q (\Omega_+)^N) \cap L_p (0, T; H^2_q (\Omega_+)^N), \\
\bu_- & \in H^1_p (0, T; L_q (\Omega_-)^N) \cap L_p (0, T; H^2_q (\Omega_-)^N), \\
\pi_- & \in L_p (0, T; H^1_q (\Omega_-)), \\
h & \in H^1_p (0, T; W^{2 - 1 \slash q}_q (S_R)) \cap L_p (0, T; W^{3 - 1 \slash q}_q (S_R)).	
\end{split} \right.	
\end{equation*}
Then, our main result in this article can be read as follows.
\begin{theo}
\label{th-main}
Let $p$ and $q$ be real numbers such that $2 < p < \infty$, $N < q < \infty$, and $2 \slash p + N \slash q < 1$. Assume that Assumptions~\ref{asp-1} and \ref{asp-2} are valid. Then, there exists a small number $\varepsilon \in (0, 1)$ such that for any initial data $\rho_{0 +} \in B^{3 - 2 \slash p}_{q, p} (\Omega_+)$, $\bu_{0 \pm} \in B^{2 (1 - 1 \slash p)}_{q, p} (\Omega_\pm)$, and $h_0 \in B^{3 - 1 \slash p - 1 \slash q}_{q, p} (S_R)$ satisfying the smallness condition:
\begin{equation*}
\lVert \rho_{0 +} \rVert_{B^{3 - 2 \slash p}_{q, p} (\Omega_+)} + \lVert \bu_{0 +} \rVert_{B^{2 (1 - 1 \slash p)}_{q, p} (\Omega_+)} + \lVert \bu_{0 -} \rVert_{B^{2 (1 - 1 \slash p)}_{q, p} (\Omega_-)} + \lVert h_0 \rVert_{B^{3 - 1 \slash p - 1 \slash q}_{q, p} (S_R)} \le \varepsilon
\end{equation*}
and the compatibility conditions:
\begin{equation}
\label{compat-1}
\left\{\begin{aligned}
\rho_{* -} \dv \bu_{0 -} = f_d (\bu_{0 -}, h_0) & = \rho_{* -} \dv \bff_d (\bu_{0 -}, h_0) & \quad & \text{in $\Omega_-$}, \\
\Pi_\bn (\mu_- \bD (\bu_-) \bn) \rvert_- - \Pi_\bn (\mu_+ \bD(\bu_+) \bn) \rvert_+ & = g (\rho_{0 +}, \bu_{0 +}, \bu_{0 -}, h_0) & \quad & \text{on $S_R$}, \\
\Pi_\bn \bu_- \rvert_- - \Pi_\bn \bu_+ \rvert_+ & = \bh (\bu_{0 +}, \bu_{0 -}, h_0) & \quad & \text{on $S_R$}, \\
\langle \nabla \rho_{0 +}, \bn \rangle \rvert_+ & = k_- (\rho_{0 +}, h_0) & \quad & \text{on $S_R$},	\\
\langle \nabla \rho_{0 +}, \bn_+ \rangle = 0, \quad \bu_{0 +} & = 0	& \quad & \text{on $\Gamma_+$}
\end{aligned}\right.
\end{equation}
and \eqref{compati-h0} the problem \eqref{eq-1.5} with $T = \infty$ admits a unique solution $(\rho_+, \bu_+, \bu_-, \pi_-, h) \in \CS_{p, q} (0, \infty)$ and the estimate
\begin{align*}
& \lVert e^{\alpha t} \pd_t \rho_+ \rVert_{L_p (0, \infty; H^1_q (\Omega_+))} + \lVert e^{\alpha t} \rho_+ \rVert_{L_p (0, \infty; H^3_q (\Omega_+))} \\
& \qquad + \sum_{\ell = \pm} \Big(\lVert e^{\alpha t} \pd_t \bu_\ell \rVert_{L_p (0, \infty; L_q (\Omega_\ell))} + \lVert e^{\alpha t} \bu_\ell \rVert_{L_p (0, \infty; H^2_q (\Omega_\ell))}\Big) \\
& \qquad + \lVert e^{\alpha t} \nabla \pi_- \rVert_{L_p (0, \infty; L_q (\Omega_-))} + \lVert e^{\alpha t} \pd_t h \rVert_{L_p (0, \infty; W^{2 - 1 \slash q}_q (S_R))} + \lVert e^{\alpha t} h \rVert_{L_p (0, \infty; W^{3 - 1 \slash q}_q (S_R))} \\
& \le C \varepsilon
\end{align*}
with some positive constants $C$ and $\alpha$ independent of $\varepsilon$.
\end{theo}

\begin{rema}
Since $\Phi (\cdot, t)$ is a $C^1$-diffeomorphism from $\Omega$ onto itself, we see that $(\varrho_+, \bv_+, \bv_-, \fp_-, \Gamma_t)$ is the unique solution to the problem \eqref{eq-1.1}--\eqref{cond-1.4} for any $t > 0$. Besides, $\varrho_+$, $\bv_\pm$, and $h$ possess the regularities
\begin{alignat*}4
\varrho_+ \in \mathrm{BUC} ([0, \infty); \mathrm{BUC}^2 (\ov{\Omega_{t +}})) \quad \bv_\pm \in \mathrm{BUC} ([0, \infty); \mathrm{BUC}^1 (\ov{\Omega_{t \pm}})), \quad h \in \mathrm{BUC} ([0, \infty); \mathrm{BUC}^2 (S_R)),
\end{alignat*}
where $\mathrm{BUC} ([0, \infty); X)$ denotes the Banach space of all $X$-valued bounded uniformly continuous functions and $\mathrm{BUC}^m (D)$ is the subset of all bounded uniformly continuous functions that has bounded partial derivatives up to order $m \in \BN$.
\end{rema}

The rest of this paper is organized as follows: In the next section, we give brief remarks on the derivation of the equations~\eqref{eq-1.5}. Section~\ref{sec-decay.estimate} is concerned with decay estimates for the linearized problem, where exponential stability of continuous analytic semigroup associated with the linearized problem is shown in Section~\ref{sec-decay}. In section~\ref{sec-nonlinear}, we prove our main result, Theorem~\ref{th-main}.

\subsection*{Notation}
Let $\BN$, $\BR$, $\BC$ be the sets of all natural numbers, real numbers, and complex numbers, respectively. Let $D \subset \BR^N$ be a domain and let $1 \le p,q \le \infty$ and $s \in \BR$. Then, $L_q (D)$, $H^{m, q} (D)$, $m \in \BN$, and $B^s_{p, q} (D)$ denote the usual Lebesgue spaces Sobolev spaces, and Besov spaces on $D$, respectively. In addition, we may write $B^s_{q, q} (D) = W^s_q (D)$ if $s \notin \BN$. For a Banach space $X$, the $m$-product space of $X$ is denoted by $X^m$ and the norm of $X^m$ is denoted by $\lVert \, \cdot \, \rVert_X$ instead of $\lVert \, \cdot \, \rVert_{X^m}$ if there is no confusion. For $I \subset \BR$ and $p \in (1, \infty]$, let $L_p (I; X)$ and $H^1_p (I; X)$ be the $X$-valued Lebesgue spaces on $I$ and the $X$-valued Sobolev spaces on $I$, respectively. Let
\begin{alignat*}4
(f, g)_D & = \int_D f (x) \ov{g (x)} \dx, & \qquad (f, g)_{S_R} & = \int_{S_R} f (x) \ov{g (x)} \,\mathrm{d} \tau, \\
(\bff, \bg)_D & = \int_D \bff (x) \cdot \ov{\bg (x)} \dx, & \qquad (\bff, \bg)_{S_R} & = \int_{S_R} \bff (x) \cdot \ov{\bg (x)} \,\mathrm{d} \tau,
\end{alignat*}
where $\mathrm{d} \tau$ is the surface element of $S_R$ and $f$ denotes the complex conjugate of $f$. The letter $C$ denotes generic constants. Besides, $C_{a, b, \dots}$ denotes a constant depending on the quantities $a$, $b$, $\dots$ The values of $C$ and $C_{a, b, \dots}$ may change from line to line.

\section{Remarks on the derivation of the equations \eqref{eq-1.5}}
\label{sec-derivation}
\noindent
Under the assumption \eqref{cond-psi}, we set
\begin{equation*}
\rho_+ (y, t) = \varrho_+ \circ \Phi - \rho_{* +}, \quad \bu_\pm (y, t) = \bv_\pm \circ \Phi, \quad \pi_- (y, t) = \fp_- \circ \Phi - \frac{\sigma (N - 1)}{R}.
\end{equation*}
Noting that the free boundary $\Gamma_t$ is given by $x = y + \Psi (y, t)$, the kinematic boundary condition reads as
\begin{equation*}
V_{\Gamma_t} = \frac{\pd x}{\pd t} \cdot \bn_t = \bigg(\frac{\pd h}{\pd t} \bn + \frac{\mathrm{d}}{\dt} \xi (t) \bigg) \cdot \bn_t.
\end{equation*}
To represent $\mathrm{d} \xi \slash \dt$, we introduce the Jacobian of $x = \Phi (y, t)$, which is denoted by $J (y, t) = 1 + J_0 (\nabla \Psi)$ with some polynomial $V (\bk)$ satisfying $V (0) = 0$. Besides, choosing $\delta$ so small, then by \eqref{cond-psi} the inverse of the Jacobi matrix of the transformation $x = \Phi (x, t)$ exists, i.e., we can write
\begin{equation*}
\bigg(\frac{\pd x}{\pd y}\bigg)^{- 1} = \bI + \sum_{k = 1}^{\infty} (- \nabla \Psi (y, t)),
\end{equation*}
and hence there exists an $N \times N$ matrix $\bV_0 (\bk)$ of $C^\infty$ functions defined on $\lvert \bk \rvert < \delta$ such that $\bV_0 (0) = 0$ and $(\pd x \slash \pd y)^{- 1} = \bI + \bV_0 (\nabla \Psi (y, t))$. Hence, it follows that
\begin{equation*}
\frac{\mathrm{d}}{\dt} \xi (t) = \frac{1}{\lvert B_R \rvert} \int_{B_R} \rho_{* -} \bu_- (y, t) \,\mathrm{d} y + \frac{1}{\lvert B_R \rvert} \int_{B_R} \rho_{* -} \bu_- (y, t) J_0 (\nabla \Psi) \,\mathrm{d} y.
\end{equation*}
Employing the similar argument given in \cite[Appendix]{W20}, we see that the functions $\rho_+$, $\bu_\pm$, $\pi_-$, and $h$ satisfy
\begin{equation}
\label{previous-fixed}
\left\{\begin{aligned}
\pd_t \rho_+ + \rho_{* +} \dv \bu_+ & = f_M (\rho_+, \bu_+, h) &\enskip &\text{ in $\Omega_+ \times(0, T)$}, \\
\rho_{* -} \dv \bu_- = f_d (\bu_-, h) & = \rho_{* -} \dv \bff_d (\bu_-, h) &\enskip &\text{ in $\Omega_- \times (0, T)$}, \\
\rho_{* +} \pd_t \bu_+ - \DV \bT_+ (\rho_+, \bu_+) & = \bff_+ (\rho_+, \bu_+, h) &\enskip &\text{ in $\Omega_+ \times (0, T)$}, \\
\rho_{* -} \pd_t \bu_- - \DV \bT_- (\bu_-, \pi_-) & = \bff_- (\bu_-, h) &\enskip &\text{ in $\Omega_- \times (0, T)$}, \\
\pd_t h - \frac{1}{\rho_{* -} - \rho_{* +}} \lbrac \langle \rho_* \bu, \bn \rangle \rbrac & = \wt d (\rho_+, \bu_+, \bu_-, h) &\enskip &\text{ on $S_R \times (0, T)$}, \\
\wt \bB (\rho_+, \bu_+, \bu_-, \pi_-, h) & = \wt \bG (\rho_+, \bu_+, \bu_-, h), &\enskip &\text{ on $S_R \times (0, T)$}, \\
\bu_+ = 0, \qquad \langle \nabla \rho_+, \bn_+ \rangle & = 0 &\enskip &\text{ on $\Gamma_+ \times (0, T)$}, \\
(\rho_+, \bu, h) \rvert_{t = 0} & = (\rho_{0 +}, \bu_0, h_0) &\enskip &\text{ on $\Omega_+ \times \dot \Omega \times S_R$},	
\end{aligned}\right.			
\end{equation}
where the right-hand members represent the nonlinear terms. See \cite[Appendix]{W20} for the exact expressions of the nonlinearities. Here, the boundary condition $\wt \bB (\rho_+, \bu_+, \bu_-, \pi_-, h) = \wt \bG (\rho_+, \bu_+, \bu_-, h)$ is given by
\begin{align}
\label{original-BC}
\left\{\begin{aligned}
\Pi_\bn (\mu_- \bD (\bu_-) \bn) \rvert_- - \Pi_\bn (\mu_+ \bD(\bu_+) \bn) \rvert_+ & = g (\rho_+, \bu_+, \bu_-, h), \\
\langle \bT_- (\bu_-, \pi_-) \bn, \bn \rangle \rvert_-
- \langle \bT_+ (\bu_+, \rho_+) \bn, \bn \rangle \rvert_+ + \sigma \CA_{S_R} h & = \wt f^+_B (\rho_+, \bu_+, \bu_-, h), \\
\frac{1}{\rho_{* -}} \langle \bT_- (\bu_-, \pi_-) \bn, \bn \rangle \rvert_- - \frac{1}{\rho_{* +}} \langle \bT_+ (\bu_+, \rho_+) \bn, \bn \rangle \rvert_+ & = \wt f^-_B (\rho_+, \bu_+, \bu_-, h),\\
\Pi_\bn \bu_- \rvert_- - \Pi_\bn \bu_+ \rvert_+ & = \bh(\bu_+, \bu_-, h), \\
\langle \nabla \rho_+, \bn \rangle \rvert_+ & = k_-(\rho_+, h),			
\end{aligned}\right.	
\end{align}

\par
We next show that the solution $(\rho_+, \bu_+, \bu_-, \pi_-, h)$ to \eqref{eq-1.5} satisfies the system \eqref{previous-fixed}. To this end, we recall Assumption~\ref{asp-1} and the representation of $\Gamma_t$. By using polar coordinates we have
\begin{align*}
\lvert B_R \rvert & = \lvert \Omega_{t -} \rvert \\
& = \int_{\lvert \omega \rvert = 1} \int_0^{R + h (R \omega, t)} r^{N - 1} \,\mathrm{d} r \,\mathrm{d} \omega \\
& = \int_{\lvert \omega \rvert = 1} \frac{1}{N} (R + h (R \omega, t))^N \,\mathrm{d} \omega, \\
& = \lvert B_R \rvert + \frac{1}{N} \sum_{j = 1}^N {}_N \mathsf{C}_j R^{1 - j} \int_{S_R} h (y, t)^j \,\mathrm{d} \tau
\intertext{and}
0 & = \frac{1}{\lvert B_R \rvert} \int_{\Omega_{t -}} (x_j - \xi_j (t)) \dx \\
& = \frac{1}{\lvert B_R \rvert} \int_{\lvert \omega \rvert = 1} \int_0^{R + h (R \omega, t)} r^N \omega_j \,\mathrm{d} r \, \mathrm{d} \omega \\
& = \int_{\lvert \omega \rvert = 1} \frac{1}{N + 1} (R + h (R \omega, t))^{N + 1} \omega_j \,\mathrm{d} \omega \\
& = (h, \varphi_j)_{S_R} + \frac{1}{N + 1} \sum_{k = 2}^{N + 1} {}_{N + 1} \mathsf{C}_k R^{1 - k} (h^k, \varphi_j)_{S_R},
\end{align*}
where $\varphi_k$ ($k = 1, \dots, N$) denote the spherical harmonics of degree 1 on $S_R$ normalized by $(\varphi_k,\varphi_\ell)_{S_R} = \delta_{k \ell}$. These formulas imply
\begin{equation}
\label{formula-h1}
\begin{split}
\int_{S_R} h (y, t) \,\mathrm{d} y & = - \frac{1}{N} \sum_{j = 2}^N {}_N \mathsf{C}_j R^{1 - j} \int_{S_R} h (y, t)^j \,\mathrm{d} \tau, \\
(h, \varphi_j)_{S_R} & = - \frac{1}{N + 1} \sum_{k = 2}^{N + 1} {}_{N + 1} \mathsf{C}_k R^{1 - k} (h^k, \varphi_j)_{S_R}.
\end{split}
\end{equation}
Especially, it holds
\begin{equation*}
\label{formula-h2}
\sum_{j = 1}^N (h, \varphi_j)_{S_R} \varphi_j = - \frac{1}{N + 1} \sum_{j = 1}^N \sum_{k = 2}^{N + 1} {}_{N + 1} \mathsf{C}_k R^{1 - k} (h^k, \varphi_j)_{S_R} \varphi_j \qquad \text{on $S_R$}.
\end{equation*}
Hence, we can rewrite \eqref{original-BC}$_{2, 3}$ as follows:
\begin{equation}
\label{boundaryboundary}
\left\{\begin{split}
\langle \bT_- (\bu_-, \pi_-) \bn, \bn \rangle \rvert_-
- \langle \bT_+ (\bu_+, \rho_+) \bn, \bn \rangle \rvert_+ + \sigma \CA_{S_R} h + \sum_{j = 1}^N (h, \varphi_j)_{S_R} \varphi_j & = f^+_B (\rho_+, \bu_+, \bu_-, h), \\
\frac{1}{\rho_{* -}} \langle \bT_- (\bu_-, \pi_-) \bn, \bn \rangle \rvert_- - \frac{1}{\rho_{* +}} \langle \bT_+ (\bu_+, \rho_+) \bn, \bn \rangle \rvert_+ + \frac{1}{\rho_{* -}} \sum_{j = 1}^N (h, \varphi_j)_{S_R} \varphi_j & = f^-_B (\rho_+, \bu_+, \bu_-, h),
\end{split}\right.
\end{equation}
where we have set
\begin{equation*}
\begin{split}
f^+_B (\rho_+, \bu_+, \bu_-, h) & = \wt f^+_B (\rho_+, \bu_+, \bu_-, h) - \frac{1}{N + 1} \sum_{j = 1}^N \sum_{k = 2}^{N + 1} {}_{N + 1} \mathsf{C}_k R^{1 - k} (h^k, \varphi_j)_{S_R} \varphi_j, \\
f^-_B (\rho_+, \bu_+, \bu_-, h) & = \wt f^-_B (\rho_+, \bu_+, \bu_-, h) - \frac{1}{(\rho_{* -} N + 1)} \sum_{j = 1}^N \sum_{k = 2}^{N + 1} {}_{N + 1} \mathsf{C}_k R^{1 - k} (h^k, \varphi_j)_{S_R} \varphi_j,
\end{split}
\end{equation*}
respectively. Notice that \eqref{boundaryboundary} is equivalent to 
\begin{equation*}
\left\{\begin{split}
\langle \bT_- (\bu_-, \pi_-) \bn, \bn \rangle \rvert_- + \frac{\rho_{* -} \sigma}{\rho_{* -} - \rho_{* +}} \CA_{S_R} h + \sum_{j = 1}^N (h, \varphi_j)_{S_R} \varphi_j & = \frac{\rho_{* -}}{\rho_{* -} - \rho_{* +}} f^+_B (\rho_+, \bu_+, \bu_-, h), \\
\langle \bT_+ (\bu_+, \rho_+) \bn, \bn \rangle \rvert_+ + \frac{\rho_{* +} \sigma}{\rho_{* -} - \rho_{* +}} \CA_{S_R} h & = \frac{\rho_{* +}}{\rho_{* -} - \rho_{* +}} f^-_B (\rho_+, \bu_+, \bu_-, h),
\end{split}\right.
\end{equation*}
\par
It remains to give a representation of the evolution equation for the height function. Let $V_{0 i j} (\bk)$ be the $(i, j)$th component of $\bV_0 (\bk)$. Then we have
\begin{equation*}
\dv \bv_+ = \sum_{j, k = 1}^N (\delta_{j k} + V_{0 j k} (\bk)) \frac{\pd u_j}{\pd y_k} = \dv \bu_+ + \bV_0 (\bk) \colon \nabla \bu_+,
\end{equation*}
and hence the transport theorem \eqref{transprot-omegaplus} can be read as
\begin{equation*}
\begin{split}
0 & = \int_{\Omega_{t +}} \dv \bv_+ \,\mathrm{d} y \\
& = \int_{\Omega_+} \Big(\dv \bu_+ + \bV_0 (\nabla \Psi) \colon \nabla \bu_+ \Big) (1 + J_0 (\nabla \Psi)) \,\mathrm{d} y.
\end{split}
\end{equation*}
Thus, the evolution equation of $h$ becomes
\begin{equation*}
\pd_t h - \frac{1}{\rho_{* -} - \rho_{* +}} \lbrac \langle \rho_* \bu, \bn \rangle \rbrac + \CM \bu = d (\rho_+, \bu_+, \bu_-, h),
\end{equation*}
where we have set
\begin{equation*}
\begin{split}
d (\rho_+, \bu_+, \bu_-, h) & = \wt d (\rho_+, \bu_+, \bu_-, h) - \frac{\rho_{* -}}{\lvert B_R \rvert} (\bn_t - \bn) \cdot \int_{B_R} \bu_- (y, t) \,\mathrm{d} y \\
& \quad + \frac{\rho_{* -}}{\lvert B_R \rvert} \bn_t \cdot \int_{B_R} \bu_- (y, t) J_0 (\nabla \Psi) \,\mathrm{d} y \\
& \quad + \frac{\rho_{* +} \varphi_0}{\rho_{* -} - \rho_{* +}} \int_{\Omega_+} \Big(\bV_0 (\nabla \Psi) \colon \nabla \bu_+ \Big) (1 + J_0 (\nabla \Psi)) \,\mathrm{d} y.
\end{split}
\end{equation*}

\section{Decay estimates for the linearized problem}
\label{sec-decay.estimate}
\noindent
To prove~\ref{th-main}, the crucial ingredient is decay properties of solutions of the Stokes equations:
\begin{align}
\label{eq-linear}
\left\{\begin{aligned}
\pd_t \rho_+ + \rho_{* +} \dv \bu_+ & = F_M &\enskip &\text{ in $\Omega_+ \times(0, T)$}, \\
\rho_{* -} \dv \bu_- = F_d & = \rho_{* -} \dv \bF_d &\enskip &\text{ in $\Omega_- \times (0, T)$}, \\
\rho_{* +} \pd_t \bu_+ - \DV \bT_+ (\rho_+, \bu_+) & = \bF_+ &\enskip &\text{ in $\Omega_+ \times (0, T)$}, \\
\rho_{* -} \pd_t \bu_- - \DV \bT_- (\bu_-, \pi_-) & = \bF_- &\enskip &\text{ in $\Omega_- \times (0, T)$}, \\
\pd_t h - \frac{1}{\rho_{* -} - \rho_{* +}} \lbrac \langle \rho_* \bu, \bn \rangle \rbrac + \CM \bu & = D &\enskip &\text{ on $S_R \times (0, T)$}, \\
\bB (\rho_+, \bu_+, \bu_-, \pi_-, h) & = \bG &\enskip &\text{ on $S_R \times (0, T)$}, \\
\bu_+ = 0, \qquad \langle \nabla \rho_+, \bn_+ \rangle & = 0 &\enskip &\text{ on $\Gamma_+ \times (0, T)$}, \\
(\rho_+, \bu, h) \rvert_{t = 0} & = (\rho_{0 +}, \bu_0, h_0) &\enskip &\text{ on $\Omega_+ \times \dot \Omega \times S_R$},
\end{aligned}\right.		
\end{align}
where $\bG = (G, F^+_B, F^-_B, \bH, K_-)$. Let $\{\varphi_\ell\}_{\ell = 0}^N$ be given by $\varphi_0 = \lvert S_R \rvert^{- 1 \slash 2}$ and $\varphi_\ell$ ($\ell = 1, \dots, N$) that denote the spherical harmonics of degree 1 on $S_R$, where $(\varphi_k,\varphi_\ell)_{S_R} = \delta_{k \ell}$. Then, $\{\varphi_\ell\}_{\ell = 0}^N$ forms an orthogonal basis of the space $\mathsf{N} (\CA_{S_R}) \cup \BC$ with respect to the $L_2 (S_R)$ inner-product $(\,\cdot \, , \, \cdot \,)_{S_R}$. In this section, we will prove the following theorem.
\begin{theo}
\label{linear-decay}
Let $1 < p, q < \infty$, $2 \slash p + 1 \slash q \ne 1$, $2 \slash p + 1 \slash q \ne 2$, and $T > 0$. There exists a constant $\eta > 0$ such that the following assertion is valid: Let $\rho_{0 +} \in B^{3 - 2 \slash p}_{q, p} (\Omega_+)$, $\bu_0 \in B^{2 (1 - 1 \slash p)}_{q, p} (\dot \Omega)$, and $h_0 \in B^{3 - 1 \slash p - 1 \slash q}_{q, p} (S_R)$. In addition, let $(F_M, F_d, \bF_d, \bF_+, \bF_-, D, G, F^+_B, F^-_B, \bH, K_-)$ be functions in the right-hand members of \eqref{eq-linear} such that
\begin{gather*}
e^{\eta t} F_M \in L_p (\BR; H^1_q (\Omega_+)), \quad e^{\eta t} F_d \in H^{1 \slash 2}_p (\BR; L_q (\Omega_-)) \cap L_p (\BR; H^1_q (\Omega_-)), \\
e^{\eta t} \bF_d \in H^1_p (\BR; L_q (\Omega_-)^N), \quad \rho_{* -} \dv \bF_d = F_d, \quad e^{\eta t} \bF_\pm \in L_p (\BR; L_q (\Omega_\pm)^N), \\
e^{\eta t} D \in L_p (\BR; W^{2 - 1 \slash q}_q (S_R)), \quad e^{\eta t} (G, F^\pm_B) \in L_p (\BR; H^1_q (\dot \Omega)) \cap H^{1 \slash 2}_p (\BR; L_q (\dot \Omega)), \\
e^{\eta t} \bH \in H^1_p (\BR; L_q (\dot \Omega)^N) \cap L_p (\BR; H^2_q (\dot \Omega)^N), \quad e^{\eta t} K_- \in L_p (\BR; H^2_q (\Omega_+)),
\end{gather*}
where the compatibility condition $\rho_{* -} \dv \bu_{0 -} = F_d \vert_{t = 0}$ holds in $\Omega_-$. Furthermore, we suppose the compatibility conditions:
\begin{align*}
\left\{\begin{aligned}
\Pi_\bn (\mu_- \bD (\bu_{0 -}) \bn) \vert_- - \Pi_\bn (\mu_+ \bD (\bu_{0 +}) \bn) \vert_+ & = G & \quad & \text{on $S_R$}, \\
\Pi_\bn \bu_{0 -} \vert_- - \Pi_\bn \bu_{0 +} \vert_+ & = \bH & \quad & \text{on $S_R$}, \\
\langle \nabla \rho_{0 +}, \bn \rangle \vert_+ & = K_- & \quad & \text{on $S_R$}, \\
\langle \nabla \rho_{0 +}, \bn_+ \rangle = 0, \quad \bu_{0 +} & = 0 & \quad & \text{on $\Gamma_+$}
\end{aligned}\right.
\end{align*}
provided $2 \slash p + 1 \slash q < 1$, while we suppose the compatibility conditions: $\langle \nabla \rho_{0 +}, \bn_+ \rangle = 0$, $\bu_{0 +} = 0$ on $\Gamma_+$ provided $2 \slash p + 1 \slash q < 2$. Then the problem \eqref{eq-linear} has a unique solution $(\rho_+, \bu_+, \bu_-, \pi_-, h) \in \CS_{p, q} (0, \infty)$ possessing the estimate
\begin{align*}
& \CJ_{p, q, T} (\rho_+, \bu_+, \bu_-, \pi_-, h; \eta) \\
& \le C \bigg\{\BJ_{p, q, T} (\bu_0, h_0, F_M, F_d, \bF_d, \bF_+, \bF_-, D, \bG; \eta) + \bigg(\int_0^T (e^{\eta s} \lvert (h (\cdot, s), \varphi_0)_{S_R} \rvert)^p \bigg)^{1 \slash p} \bigg\}
\end{align*}
for some constant independent of $\eta$ and $T$. Here and in the following, we set
\begin{align*}
& \CJ_{p, q, T} (\rho_+, \bu_+, \bu_-, \pi_-, h; \eta) \\
& \quad \quad = \lVert e^{\eta t} \pd_t \rho_+ \rVert_{L_p (0, T; H^1_q (\Omega_+))} + \lVert e^{\eta t} \rho_+ \rVert_{L_p (0, T; H^3_q (\Omega_+))} \\
& \quad \quad \quad + \sum_{\ell = \pm} \Big(\lVert e^{\eta t} \pd_t \bu_\ell \rVert_{L_p (0, T; L_q (\Omega_\ell))} + \lVert e^{\eta t} \bu_\ell \rVert_{L_p (0, T; H^2_q (\Omega_\ell))}\Big) \\
& \quad \quad \quad + \lVert e^{\eta t} \nabla \pi_- \rVert_{L_p (0, T; L_q (\Omega_-))} + \lVert e^{\eta t} \pd_t h \rVert_{L_p (0, T; W^{2 - 1 \slash q}_q (S_R))} + \lVert e^{\eta t} h \rVert_{L_p (0, T; W^{3 - 1 \slash q}_q (S_R))}, \\
& \BJ_{p, q, T} (\bu_0, h_0, F_M, F_d, \bF_d, \bF_+, \bF_-, D, \bG; \eta) \\
& \quad \quad = \lVert \rho_{0 +} \rVert_{B^{3 - 2 \slash p}_{q, p} (\Omega_+)} + \sum_{\ell = \pm} \lVert \bu_{0 \ell} \rVert_{B^{2 (1 - 1 \slash p)}_{q, p} (\Omega_\ell)} + \lVert h_0 \rVert_{B^{3 - 1 \slash p - 1 \slash q}_{q, p} (S_R)} + \lVert e^{\eta t} F_M \rVert_{L_p (\BR; H^1_q (\Omega_+))} \\
& \quad \quad \quad + \lVert e^{\eta t} F_d \rVert_{L_p (\BR; H^1_q (\Omega_-))} + \lVert e^{\eta t} F_d \rVert_{H^{1 \slash 2}_p (\BR; L_q (\Omega_-))} + \lVert e^{\eta t} \pd_t \bF_d \rVert_{L_p (\BR; L_q (\Omega_-))} + \sum_{\ell = \pm} \lVert e^{\eta t} \bF_\ell \rVert_{L_p (\BR; L_q (\Omega_\ell))} \\
& \quad \quad \quad + \lVert e^{\eta t} D \rVert_{L_p (\BR; W^{2 - 1 \slash q}_q (S_R))} + \lVert e^{\eta t} (G, \nabla \bH) \rVert_{H^{1 \slash 2}_p (\BR; L_q (\dot \Omega))} + \lVert e^{\eta t} (\nabla G, \pd_t \bH, \nabla \bH) \rVert_{L_p (\BR; L_q (\dot \Omega))} \\
& \quad \quad \quad + \sum_{\ell = \pm} \Big(\lVert e^{\eta t} F^\ell_B \rVert_{L_p (\BR; L_q (\dot \Omega))} + \lVert e^{\eta t} F^\ell_B \rVert_{H^{1 \slash 2}_p (\BR; L_q (\dot \Omega))} + \lVert e^{\eta t} \nabla F^\ell_B \rVert_{L_p (\BR; L_q (\dot \Omega))} \Big) \\
& \quad \quad \quad + \lVert e^{\eta t} \nabla K_- \rVert_{H^{1 \slash 2}_p (\BR; L_q (\dot \Omega))} + \lVert e^{\eta t} (\pd_t K_-, \nabla^2 K_-) \rVert_{L_p (\BR; L_q (\dot \Omega))}.
\end{align*}
\end{theo}
To show Theorem~\ref{linear-decay}, we first consider the shifted equations
\begin{align}
\label{eq-shifted}
\left\{\begin{aligned}
\pd_t \rho_{1 +} + \lambda_1 \rho_{1 +} + \rho_{* +} \dv \bu_{1 +} & = F_M &\enskip &\text{ in $\Omega_+ \times(0, T)$}, \\
\rho_{* -} \dv \bu_{1 -} = F_d & = \rho_{* -} \dv \bF_d &\enskip &\text{ in $\Omega_- \times (0, T)$}, \\
\rho_{* +} \pd_t \bu_{1 +} + \rho_{* +} \lambda_1 \bu_{1 +} - \DV \bT_+ (\rho_{1 +}, \bu_{1 +}) & = \bF_+ &\enskip &\text{ in $\Omega_+ \times (0, T)$}, \\
\rho_{* -} \pd_t \bu_{1 -} + \rho_{* +} \lambda_1 \bu_{1 -} - \DV \bT_- (\bu_{1 -}, \pi_{1 -}) & = \bF_- &\enskip &\text{ in $\Omega_- \times (0, T)$}, \\
\pd_t h_1 + \lambda_1 h - \frac{1}{\rho_{* -} - \rho_{* +}} \lbrac \langle \rho_* \bu_1, \bn \rangle \rbrac + \CM \bu_1 & = D &\enskip &\text{ on $S_R \times (0, T)$}, \\
\bB (\rho_{1 +}, \bu_{1 +}, \bu_{1 -}, \pi_{1 -}, h_1) & = \bG &\enskip &\text{ on $S_R \times (0, T)$}, \\
\bu_{1 +} = 0, \qquad \langle \nabla \rho_{1 +}, \bn_+ \rangle & = 0 &\enskip &\text{ on $\Gamma_+ \times (0, T)$}, \\
(\rho_{1 +}, \bu_1, h_1) \rvert_{t = 0} & = (\rho_{0 +}, \bu_0, h_0) &\enskip &\text{ on $\Omega_+ \times \dot \Omega \times S_R$}.	
\end{aligned}\right.			
\end{align}
For the shifted equation~\eqref{eq-shifted}, the following theorem can be shown.
\begin{theo}
\label{th-shifted}
Let $1 < p, q < \infty$, $2 \slash p + 1 \slash q \ne 1$, $2 \slash p + 1 \slash q \ne 2$, and $T > 0$. Then, there exists a constant $\lambda_1 > 0$ such that if $\lambda_1 \ge \lambda_2$, then the following assertion holds: Let $\rho_{0 +} \in B^{3 - 2 \slash p}_{q, p} (\Omega_+)$, $\bu_0 \in B^{2 (1 - 1 \slash p)}_{q, p} (\dot \Omega)^N$, and $h_0 \in B^{3 - 1 \slash p - 1 \slash q}_{q, p} (S_R)$ be the initial data for equations~\eqref{eq-shifted} and let $F_M$, $F_d$, $\bF_d$, $\bF_\pm$, $D$, $\bG$ be given functions on the right-hand side of~\eqref{eq-shifted} with
\begin{gather*}
F_M \in L_p (\BR; H^1_q (\Omega_+)), \quad F_d \in H^{1 \slash 2}_p (\BR; L_q (\Omega_-)) \cap L_p (\BR; H^1_q (\Omega_-)), \\
\bF_d \in H^1_p (\BR; L_q (\Omega_-)^N), \quad \rho_{* -} \dv \bF_d = F_d, \quad  \bF_\pm \in L_p (\BR; L_q (\Omega_\pm)^N), \\
D \in L_p (\BR; W^{2 - 1 \slash q}_q (S_R)), \quad G, F^\pm_B \in H^{1 \slash 2}_p (\BR; L_q (\dot \Omega)) \cap L_p (\BR; H^1_q (\dot \Omega)), \\
\bH \in H^1_p (\BR; L_q (\dot \Omega)^N) \cap L_p (\BR; H^2_q (\dot \Omega)^N), \quad K_- \in L_p (\BR; H^2_q (\Omega_+)),	
\end{gather*}
where the compatibility condition $\rho_{* -} \dv \bu_{0 -} = F_d \vert_{t = 0}$ is valid in $\Omega_-$. Suppose the compatibility conditions:
\begin{align*}
\left\{\begin{aligned}
\Pi_\bn (\mu_- \bD (\bu_{0 -}) \bn) \vert_- - \Pi_\bn (\mu_+ \bD (\bu_{0 +}) \bn) \vert_+ & = G & \quad & \text{on $S_R$}, \\
\Pi_\bn \bu_{0 -} \vert_- - \Pi_\bn \bu_{0 +} \vert_+ & = \bH & \quad & \text{on $S_R$}, \\
\langle \nabla \rho_{0 +}, \bn \rangle \vert_+ & = K_- & \quad & \text{on $S_R$}, \\
\langle \nabla \rho_{0 +}, \bn_+ \rangle = 0, \quad \bu_{0 +} & = 0 & \quad & \text{on $\Gamma_+$}		
\end{aligned}\right.
\end{align*}
if $2 \slash p + 1 \slash q < 1$, while
\begin{align*}
\langle \nabla \rho_{0 +}, \bn_+ \rangle = 0, \quad \bu_{0 +} = 0 \quad \text{on $\Gamma_+$}	
\end{align*}
if $2 \slash p + 1 \slash q < 2$. Then, the problem~\eqref{eq-shifted} admits a unique solution $(\rho_{1 +}, \bu_{1 +}, \bu_{1 -}, \pi_{1 -}, h_1) \in \CS_{p, q} (0, \infty)$ possessing the estimate
\begin{equation}
\label{est-shift}
\CJ_{p, q, T} (\rho_{1 +}, \bu_{1 +}, \bu_{1 -}, \pi_{1 -}, h_1; 0) \le C \BJ_{p, q, T} (\bu_0, h_0, F_M, F_d, \bF_d, \bF_+, \bF_-, D, \bG; 0)
\end{equation}
for some constant $C$ independent of $T$.
\end{theo}
\begin{proof}
Employing the argument in \cite{W20}, we can show the unique existence of $(\rho_{1 +}, \bu_{1 +}, \bu_{1 -}, \pi_{1 -}, h_1)$ possessing the estimate \eqref{est-shift}. In fact, we can show the existence of the $\CR$-bounded solution operators for the generalized resolvent problem that is obtained by the Laplace transform of \eqref{eq-shifted} with respect to time $t$, and hence the operator valued Fourier multiplier theorem (cf. Weis~\cite{W01}) yields the estimate \eqref{est-shift}. We refer to \cite{W20} for the detailed proof.
\end{proof}
For any $\eta > 0$, we see that $(e^{\eta t} \rho_{1 +}, e^{\eta t} \bu_{1 +}, e^{\eta t} \bu_{1 -}, e^{\eta t} \pi_{1 -}, e^{\eta t} h_1)$ satisfies the equations:
\begin{align*}
\left\{\begin{aligned}
\pd_t (e^{\eta t} \rho_{1 +}) + (\lambda_1 - \eta) (e^{\eta t} \rho_{1 +}) + \rho_{* +} \dv (e^{\eta t} \bu_{1 +})& = e^{\eta t} F_M &\enskip &\text{ in $\Omega_+ \times(0, T)$}, \\
\rho_{* -} \dv (e^{\eta t} \bu_{1 -}) = e^{\eta t} F_d & = \rho_{* -} \dv e^{\eta t} \bF_d &\enskip &\text{ in $\Omega_- \times (0, T)$}, \\
\rho_{* +} \pd_t (e^{\eta t} \bu_{1 +}) + \rho_{* +} (\lambda_1 - \eta) (e^{\eta t} \bu_{1 +}) - \DV \bT_+ (e^{\eta t} \rho_{1 +}, e^{\eta t} \bu_{1 +}) & = e^{\eta t} \bF_+ &\enskip &\text{ in $\Omega_+ \times (0, T)$}, \\
\rho_{* -} \pd_t (e^{\eta t} \bu_{1 -}) + \rho_{* +} (\lambda_1 - \eta) e^{\eta t} \bu_{1 -} - \DV \bT_- (\bu_{1 -}, \pi_{1 -}) & = e^{\eta t} \bF_- &\enskip &\text{ in $\Omega_- \times (0, T)$}, \\
\pd_t (e^{\eta t} h_1) + (\lambda_1 - \eta) (e^{\eta t} h_1) \qquad \\
- \frac{1}{\rho_{* -} - \rho_{* +}} \lbrac \langle \rho_* e^{\eta t} \bu_1, \bn \rangle \rbrac + \CM (e^{\eta t} \bu_1) & = e^{\eta t} D &\enskip &\text{ on $S_R \times (0, T)$}, \\
\bB (e^{\eta t} \rho_{1 +}, e^{\eta t} \bu_{1 +}, e^{\eta t} \bu_{1 -}, e^{\eta t} \pi_{1 -}, e^{\eta t} h_1) & = e^{\eta t} \bG &\enskip &\text{ on $S_R \times (0, T)$}, \\
e^{\eta t} \bu_{1 +} = 0, \qquad \langle \nabla (e^{\eta t} \rho_{1 +}), \bn_+ \rangle & = 0 &\enskip &\text{ on $\Gamma_+ \times (0, T)$}, \\
(e^{\eta t} \rho_{1 +}, e^{\eta t} \bu_1, e^{\eta t} h_1) \rvert_{t = 0} & = (\rho_{0 +}, \bu_0, h_0) &\enskip &\text{ on $\Omega_+ \times \dot \Omega \times S_R$}.		
\end{aligned}\right.				
\end{align*}
For given $\eta > 0$ we choose $\lambda_1 > 0 $ such that $\lambda_1 - \eta \ge \lambda_2$, from Theorem~\ref{th-shifted}, we obtain the next corollary.
\begin{corr}
\label{cor-shifted}
Let $1 < p, q < \infty$, $2 \slash p + 1 \slash q \ne 1$, $2 \slash p + 1 \slash q \ne 2$, $\eta > 0$, and $T > 0$. Let $\rho_{0 +} \in B^{3 - 2 \slash p}_{q, p} (\Omega_+)$, $\bu_0 \in B^{2 (1 - 1 \slash p)}_{q, p} (\dot \Omega)$, and $h_0 \in B^{3 - 1 \slash p - 1 \slash q}_{q, p} (S_R)$. In addition, let $(F_M, F_d, \bF_d, \bF_+, \bF_-, D, G, F^+_B, F^-_B, \bH, K_-)$ be functions in the right-hand members of \eqref{eq-shifted} such that
{\allowdisplaybreaks
\begin{gather*}
e^{\eta t} F_M \in L_p (\BR; H^1_q (\Omega_+)), \quad e^{\eta t} F_d \in H^{1 \slash 2}_p (\BR; L_q (\Omega_-)) \cap L_p (\BR; H^1_q (\Omega_-)), \\
e^{\eta t} \bF_d \in H^1_p (\BR; L_q (\Omega_-)^N), \quad \rho_{* -} \dv \bF_d = F_d, \quad e^{\eta t} \bF_\pm \in L_p (\BR; L_q (\Omega_\pm)^N), \\
e^{\eta t} D \in L_p (\BR; W^{2 - 1 \slash q}_q (S_R)), \quad e^{\eta t} (G, F^\pm_B) \in  H^{1 \slash 2}_p (\BR; L_q (\dot \Omega)) \cap L_p (\BR; H^1_q (\dot \Omega)), \\
e^{\eta t} \bH \in H^1_p (\BR; L_q (\dot \Omega)^N) \cap L_p (\BR; H^2_q (\dot \Omega)^N), \quad e^{\eta t} K_- \in L_p (\BR; H^2_q (\Omega_+)),	
\end{gather*}} where the compatibility condition $\rho_{* -} \dv \bu_{0 -} = F_d \vert_{t = 0}$ holds in $\Omega_-$. Furthermore, we suppose the compatibility conditions:
\begin{align*}
\left\{\begin{aligned}
\Pi_\bn (\mu_- \bD (\bu_{0 -}) \bn) \vert_- - \Pi_\bn (\mu_+ \bD (\bu_{0 +}) \bn) \vert_+ & = G & \quad & \text{on $S_R$}, \\
\Pi_\bn \bu_{0 -} \vert_- - \Pi_\bn \bu_{0 +} \vert_+ & = \bH & \quad & \text{on $S_R$}, \\
\langle \nabla \rho_{0 +}, \bn \rangle \vert_+ & = K_- & \quad & \text{on $S_R$}, \\
\langle \nabla \rho_{0 +}, \bn_+ \rangle = 0, \quad \bu_{0 +} & = 0 & \quad & \text{on $\Gamma_+$}	
\end{aligned}\right.	
\end{align*}
provided $2 \slash p + 1 \slash q < 1$, while we suppose the compatibility conditions:
\begin{align*}
\langle \nabla \rho_{0 +}, \bn_+ \rangle = 0, \quad \bu_{0 +} = 0 \quad \text{on $\Gamma_+$}	
\end{align*}
provided $2 \slash p + 1 \slash q < 2$. Then, there exists a constant $\lambda_1 > 0$ such that the problem \eqref{eq-shifted} admits a unique solution $(\rho_{1 +}, \bu_{1 +}, \bu_{1 -}, \pi_{1 -}, h_1) \in \CS_{p, q} (0, \infty)$ possessing the estimate
\begin{equation}
\label{shift-decay}
\CJ_{p, q, T} (\rho_{1 +}, \bu_{1 +}, \bu_{1 -}, \pi_{1 -}, h_1; \eta) \le C \BJ_{p, q, T} (\bu_0, h_0, F_M, F_d, \bF_d, \bF_+, \bF_-, D, \bG; \eta)
\end{equation}
for some constant independent of $\eta$ and $T$.	
\end{corr}
We seek the solution $(\rho_+, \bu_+, \bu_-, \pi_-, h)$ of \eqref{eq-linear} of the form: $\rho_+ = \rho_{1 +} + \rho_{2 +}$, $\bu_\pm = \bu_{1 \pm} + \bu_{2 \pm}$, $\pi_- = \pi_{1 -} + \pi_{2 -}$, and $h = h_1 + h_2$, where $(\rho_{1 +}, \bu_{1 +}, \bu_{1 -}, \pi_{1 -}, h_1)$ enjoys the shifted equation~\eqref{eq-shifted}. Then, we find that $(\rho_{2 +}, \bu_{2 +}, \bu_{2 -}, \pi_{2 -}, h_2)$ satisfies
\begin{equation}
\label{eq-homo}
\left\{\begin{aligned}
\pd_t \rho_{2 +} + \rho_{* +} \dv \bu_{2 +} & = - \lambda_1 \rho_{1 +} &\enskip &\text{ in $\Omega_+ \times(0, T)$}, \\
\rho_{* -} \dv \bu_{2 -} & = 0 &\enskip &\text{ in $\Omega_- \times (0, T)$}, \\
\rho_{* +} \pd_t \bu_{2 +} - \DV \bT_+ (\rho_{2 +}, \bu_{2 +}) & = - \rho_{* +} \lambda_1 \bu_{1 +} &\enskip &\text{ in $\Omega_+ \times (0, T)$}, \\
\rho_{* -} \pd_t \bu_{2 -} - \DV \bT_- (\bu_{2 -}, \pi_{2 -}) & = - \rho_{* +} \lambda_1 \bu_{1 -} &\enskip &\text{ in $\Omega_- \times (0, T)$}, \\
\pd_t h_2 - \frac{1}{\rho_{* -} - \rho_{* +}} \lbrac \langle \rho_* \bu_2, \bn \rangle \rbrac + \CM \bu_2 & = - \lambda_1 h_1 &\enskip &\text{ on $S_R \times (0, T)$}, \\
\bB (\rho_{2 +}, \bu_{2 +}, \bu_{2 -}, \pi_{2 -}, h_2) & = 0 &\enskip &\text{ on $S_R \times (0, T)$}, \\
\bu_{2 +} = 0, \qquad \langle \nabla \rho_{2 +}, \bn_+ \rangle & = 0 &\enskip &\text{ on $\Gamma_+ \times (0, T)$}, \\
(\rho_{2 +}, \bu_2, h_2) \rvert_{t = 0} & = (0, 0, 0) &\enskip &\text{ on $\Omega_+ \times \dot \Omega \times S_R$}.		
\end{aligned}\right.	
\end{equation}
\par
For $1 < q < \infty$, we define solenoidal spaces $J_q (\Omega_-)$ by
\begin{equation*}
J_q (\Omega_-) = \{\bF_- \in L_q (\Omega_-)^N \,\colon\, (\bF_-, \nabla \varphi)_{\Omega_-} = 0 \enskip \text{for any $\varphi \in \wh H^1_{q', 0} (\Omega_-)$}\},
\end{equation*}
where $\wh H^1_{q', 0} (\Omega_-)$ is defined by $\wh H^1_{q', 0} (\Omega_-) = \{\theta \in L_{q', \mathrm{loc}} (\Omega_-) \,\colon\, \nabla \theta \in L_{q'} (\Omega_-)^N, \enskip \theta \vert_{S_R} = 0\}$ with $q' = q \slash (q - 1)$. Since $\Omega_-$ is bounded, we know that $C^\infty_0 (\Omega_-)$ is dense in $\wh H^1_{q', 0} (\Omega_-)$, so that the necessary and sufficient condition in order that $\bu_{2 -} \in J_q (\Omega_-)$ is that $\dv \bu_{2 -} = 0$ in $\Omega_-$. Define
\begin{equation*}
H^1_{q, 0} (\Omega_-) = \{\theta \in H^1_q (\Omega_-) \,\colon\, \theta \vert_{S_R} = 0\}
\end{equation*}
and let $\psi \in H^1_{q, 0} (\Omega_-)$ be a solution to the variational equation
\begin{equation*}
(\nabla \psi, \nabla \varphi)_{\Omega_-} = (\bu_{1 -}, \nabla \varphi)_{\Omega_-} \quad \text{for any $\varphi \in H^1_{q', 0} (\Omega_-)$}.
\end{equation*}
Set $\bw = \bu_{1 -} - \nabla \psi$. Then, it holds $\bw \in J_q (\Omega_-)$ and
\begin{equation*}
\lVert \bw \rVert_{L_q (\Omega_-)} + \lVert \nabla \psi \rVert_{L_q (\Omega_-)} \le C \lVert \bu_{1 -} \rVert_{L_q (\Omega_-)}.
\end{equation*}
Using $\bw$ and $\psi$, the equations~\eqref{eq-homo}$_{2, 4}$ can be rewritten as
\begin{equation*}
\rho_{* -} \pd_t \bu_{2 -} - \DV \bT_- (\bu_{2 -}, \pi_{2 -} + \rho_{* -} \lambda_1 \nabla \psi) = - \lambda_1 \bw, \quad \dv \bu_{2 -} = 0 \qquad \text{in $\Omega_- \times (0, T)$}.
\end{equation*}
From this viewpoint, in what follows, we may suppose that
\begin{equation*}
\bu_{1 -} \in H^1_p (0, T; J_q (\Omega_-)) \cap L_p (0, T; H^2_q (\Omega_-)^N).
\end{equation*}
\par
Notice that the second and third boundary condition of $\bB (\rho_{2 +}, \bu_{2 +}, \bu_{2 -}, \pi_{2 -}, h_2) = 0$ are equivalent to
\begin{equation*}
\left\{\begin{split}
\langle \bT_- (\bu_{2 -}, \pi_{2 -}) \bn, \bn \rangle \bigg\vert_- + \frac{\rho_{* -} \sigma}{\rho_{* -} + \rho_{* +}} \CA_{S_R} h_2 + \sum_{j = 1}^N (h_2, \varphi_j)_{S_R} \varphi_j & = 0 \qquad \text{on $S_R$}, \\
\langle \bT_+ (\rho_{2 +}, \bu_{2 +}) \bn, \bn \rangle \bigg\vert_+ + \frac{\rho_{* +} \sigma}{\rho_{* -} - \rho_{* +}} \CA_{S_R} h_2 & = 0 \qquad \text{on $S_R$}.	
\end{split}\right.
\end{equation*}
Recalling the argument in Sections~4.1 in~\cite{W20}, we introduce a functional $\mathsf{P} (\bu_{2 -}, h_2) \in H^1_q (\Omega_-) + \wh H^1_q (\Omega_-)$ that is a unique solution to the weak problem
\begin{equation}
\label{pressure}
(\nabla \mathsf{P} (\bu_{2 -}, h), \nabla \varphi)_{\Omega_-} = (\DV (\mu_- \bD (\bu_{2 -})) - \rho_{* -} \nabla \dv \bu_{2 -}, \nabla \varphi)_{\Omega_-}
\end{equation}
for any $\varphi \in \wh H^1_{q', 0} (\Omega_-)$ subject to
\begin{align*}
\mathsf{P} (\bu_{2 -}, h_2) & = \mu_- \langle \bD (\bu_{2 -}) \bn, \bn \rangle + \frac{\rho_{* -} \sigma}{\rho_{* -} - \rho_{* +}} \CA_{S_R} h_2 \\
& \quad + \sum_{j = 1}^N (h_2, \varphi_j)_{S_R} \varphi_j - \rho_{* -} \dv \bu_{2 -} \qquad\qquad\qquad \text{on $S_R$}.
\end{align*}
To handle \eqref{eq-homo} in the semigroup setting, we consider the initial value problem:
\begin{equation}
\label{eq-initial}
\left\{\begin{aligned}
\pd_t \rho_{2 +} + \rho_{* +} \dv \bu_{2 +} & = 0 &\enskip &\text{ in $\Omega_+ \times(0, T)$}, \\
\rho_{* -} \dv \bu_{2 -} & = 0 &\enskip &\text{ in $\Omega_- \times (0, T)$}, \\
\rho_{* +} \pd_t \bu_{2 +} - \DV \bT_+ (\rho_{2 +}, \bu_{2 +}) & = 0 &\enskip &\text{ in $\Omega_+ \times (0, T)$}, \\
\rho_{* -} \pd_t \bu_{2 -} - \DV \bT_- (\bu_{2 -}, \mathsf{P} (\bu_{2 -}, h_2)) & = 0 &\enskip &\text{ in $\Omega_- \times (0, T)$}, \\
\pd_t h_2 - \frac{1}{\rho_{* -} - \rho_{* +}} \lbrac \langle \rho_* \bu_2, \bn \rangle \rbrac + \CM \bu_2 & = 0 &\enskip &\text{ on $S_R \times (0, T)$}, \\
\bB (\rho_{2 +}, \bu_{2 +}, \bu_{2 -}, \mathsf{P} (\bu_{2 -}, h_2), h_2) & = 0 &\enskip &\text{ on $S_R \times (0, T)$}, \\
\bu_{2 +} = 0, \qquad \langle \nabla \rho_{2 +}, \bn_+ \rangle & = 0 &\enskip &\text{ on $\Gamma_+ \times (0, T)$}, \\
(\rho_{2 +}, \bu_2, h_2) \rvert_{t = 0} & = (\rho_{2 0 +}, \bu_{2 0}, h_{2 0}) &\enskip &\text{ on $\Omega_+ \times \dot \Omega \times S_R$}.			
\end{aligned}\right.	
\end{equation}
Here, from the definition of $\mathsf{P} (\bu_{2 -}, h_2)$, we observe that the boundary condition \eqref{eq-initial}$_6$ can be written as
\begin{align}
\label{initial-BC}
\left\{\begin{aligned}
\Pi_\bn (\mu_- \bD (\bu_{2 -}) \bn) \rvert_- - \Pi_\bn (\mu_+ \bD (\bu_{2 +}) \bn) \rvert_+ & = 0, \\
\rho_{* -} \dv \bu_{2 -} & = 0, \\
\langle \bT_+ (\rho_{2 +}, \bu_{2 +}) \bn, \bn \rangle \bigg\vert_+ - \frac{\rho_{* +} \sigma}{\rho_{* -} - \rho_{* +}} \CA_{S_R} h_2 & = 0,\\
\Pi_\bn \bu_{2 -} \rvert_- - \Pi_\bn \bu_{2 +} \rvert_+ & = 0, \\
\langle \nabla \rho_{2 +}, \bn \rangle \rvert_+ & = 0.			
\end{aligned}\right.	
\end{align}
Now, we define $\mathsf{B}_q$, $\mathsf{D}_q$, and $A_q$ by
\begin{equation*}
\begin{split}
\mathsf{B}_q & = \{(\rho_{2 +}, \bu_{2 +}, \bu_{2 -}, h_2) \,\colon\, \rho_{2 +} \in H^1_q (\Omega_+), \enskip \bu_{2 +} \in L_q (\Omega_+)^N, \enskip \bu_{2 -} \in J_q (\Omega_-), \enskip h_2 \in W^{2 -  1 \slash q}_q (S_R)\}, \\
\mathsf{D}_q & = \left\{(\rho_{2 +}, \bu_{2 +}, \bu_{2 -}, h_2) \,\colon\, \begin{array}{l}
\rho_{2 +} \in H^3_q (\Omega_+), \enskip \bu_{2 +} \in H^2_q (\Omega_+)^N, \enskip \bu_{2 -} \in H^2_q (\Omega_-)^N, \\
h_2 \in W^{3 -  1 \slash q}_q (S_R), \enskip \text{$(\rho_{2 +}, \bu_{2 +}, \bu_{2 -}, h_2)$ satisfies \eqref{initial-BC} on $S_R \times (0, T)$} \\
\text{and the Dirichlet boundary conditons: $\bu_{2 +} = 0$ and $\langle \rho_{2 +}, \bn_+ \rangle = 0$} \\
\text{on $\Gamma_+ \times (0, T)$}
\end{array}\right\},
\end{split}	
\end{equation*}
and
\begin{equation*}
\label{operator-A}
A_q (\rho_{2 +}, \bu_2, h_2) = \begin{pmatrix}
- \rho_{* +} \dv \bu_{2 +} \\
\rho_{* +}^{- 1} \DV \bT_+ (\rho_{2 +}, \bu_{2 +}) \\
\rho_{* -}^{- 1} \DV \bT_- (\bu_{2 -}, \mathsf{P} (\bu_{2 -}, h_2)) \\
(\rho_{* -} - \rho_{* +})^{- 1} \lbrac \langle \rho_* \bu_2, \bn \rangle \rbrac - \CM \bu_2
\end{pmatrix},
\end{equation*}
respectively.
Then, \eqref{eq-initial} is formulated by
\begin{equation}
\label{Cauchy}
\begin{split}
\pd_t \, \CU (t) - A_q \, \CU (t) & = 0 \quad (t > 0), \\
\CU (0) & = \CU_0
\end{split}
\end{equation}
with $\CU = (\rho_{2 +}, \bu_2, h_2) \in \mathsf{D}_q$ for $t > 0$ and $\CU_0 = (\rho_{2 0 +}, \bu_{2 0}, h_{20}) \in \mathsf{B}_q$. According to \cite[Thm.~6.3]{W20}, the operator $A_q$ generates a $C^0$-analytic semigroup $\{e^{- A_q t}\}_{t \ge 0}$ on $\mathsf{B}_q$. As usual, the resolvent of $A_q$ is denoted by $R (\lambda, A_q)$. Furthermore, for some $\lambda_0$ and $0 < \varepsilon_0 < \pi \slash 2$, the set
\begin{equation}
\label{sector}
\Sigma_{\varepsilon_0, \lambda_0} := \{\lambda \in \BC \,\colon\, \lvert \arg (\lambda) \rvert \le \pi - \varepsilon_0, \enskip \lvert \lambda \rvert \ge \lambda_0\}
\end{equation}
is included in the resolvent set of $A_q$. We define
\begin{equation*}
\dot{\mathsf{B}}_q = \{(\rho_{2 +}, \bu_{2 +}, \bu_{2 -}, h_2) \in \mathsf{B}_q \,\colon\,  \enskip (h_2, 1)_{S_R} = 0\}.
\end{equation*}
To address some exponential decay property of \eqref{eq-initial}, we here record some simple but important fact that the closed subspace $\dot{\mathsf{B}}_q$ is $e^{- A_q t}$-invariant, i.e., $e^{- A_q t} \dot{\mathsf{B}}_q \subset \dot{\mathsf{B}}_q$ for any $t \ge 0$.
\begin{lemm}
\label{lem-invariant}
Let $\Sigma_{\varepsilon_0, \lambda_0}$ be given in \eqref{sector}. Then the subspace $\dot{\mathsf{B}}_q$ is $R (\lambda, A_q)$-invariant for every $\lambda \in \Sigma_{\varepsilon_0, \lambda_0}$. Namely, for given $(F_M, \bF_+, \bF_-, D) \in \dot{\mathsf{B}}_q$, the solution $(\rho_+, \bu_+, \bu_-, h)$ to  
\begin{equation}
\label{resolvent-invariant}
\left\{\begin{aligned}
\lambda \rho_+ + \rho_{* +} \dv \bu_+ & = F_M &\enskip &\text{ in $\Omega_+$}, \\
\rho_{* -} \dv \bu_- & = 0 &\enskip &\text{ in $\Omega_-$}, \\
\rho_{* +} \lambda \bu_+ - \DV \bT_+ (\rho_+, \bu_+) & = \bF_+ &\enskip &\text{ in $\Omega_+$}, \\
\rho_{* -} \lambda \bu_- - \DV \bT_- (\bu_-, \mathsf{P} (\bu_-, h)) & = \bF_- &\enskip &\text{ in $\Omega_-$}, \\
\lambda h - \frac{1}{\rho_{* -} - \rho_{* +}} \lbrac \langle \rho_* \bu, \bn \rangle \rbrac + \CM \bu & = D &\enskip &\text{ on $S_R$}, \\
\bB (\rho_+, \bu_+, \bu_-, \mathsf{P} (\bu_-, h)) & = 0 &\enskip &\text{ on $S_R$}, \\
\bu_+ = 0, \qquad \langle \nabla \rho_+, \bn_+ \rangle & = 0 &\enskip &\text{ on $\Gamma_+$}.				
\end{aligned}\right.			
\end{equation}
belongs to $\dot{\mathsf{B}}_q$.
\end{lemm}
\begin{proof}
It suffices to prove that $(h, 1)_{S_R} = 0$. Integrating \eqref{resolvent-invariant}$_4$, we have $\lambda (h, 1)_{S_R} = 0$ since $\dv \bu_- = 0$. Noting that $\lambda \ne 0$, we obtain $(h, 1)_{S_R} = 0$. This completes the proof.
\end{proof}
According to Lemma~\ref{lem-invariant}, the restriction operator $\wt A_q := A_q \vert_{\dot{\mathsf{B}}_q}$ with its domain given by $\mathsf{D} (\wt A_q) = \mathsf{D} (A_q) \cap \dot{\mathsf{B}}_q$ is the generator of the induced $C^0$-semigroup $\{e^{- \wt A_q t}\}_{t \ge 0} := \{e^{- A_q t} \vert_{\dot{\mathsf{B}}_q}\}_{t \ge 0}$, which is analytic. Since $\Omega$ is bounded, we can show that $\lambda = 0$ is included in the resolvent set of $\wt A_q$, which implies that the induced $C^0$-semigroup $\{e^{- \wt A_q t}\}_{t \ge 0}$ is exponentially stable on $\dot{\mathsf{B}}_q$.
\begin{theo}
\label{th-semigroup}
Let $1 < q < \infty$. Then, the induced $C^0$-semigroup $\{e^{- \wt A_q t}\}_{t \ge 0}$ is exponentially stable on $\dot{\mathsf{B}}_q$, that is,
\begin{equation*}
\lVert e^{- \wt A_q t} (f_1, \bff_2, \bff_3, f_4) \rVert_{\mathsf{B}_q} \le C e^{- \eta_1 t} \lVert (f_1, \bff_2, \bff_3, f_4) \rVert_{\mathsf{B}_q}
\end{equation*}
for any $t > 0$ and $(f_1, \bff_2, \bff_3, f_4) \in \dot{\mathsf{B}}_q$ with some positive constants $C$ and $\eta_1$. Here, we have set
\begin{equation*}
\lVert (f_1, \bff_2, \bff_3, f_4) \rVert_{\mathsf{B}_q} = \lVert f_1 \rVert_{H^1_q (\Omega_+)} + \lVert \bff_2 \rVert_{L_q (\Omega_+)} + \lVert \bff_3 \rVert_{L_q (\Omega_-)} + \lVert f_4 \rVert_{W^{2 - 1 \slash q}_q (S_R)}.
\end{equation*}	
\end{theo}
We will give the proof of Theorem~\ref{th-semigroup} in the next section and we now continue the proof of Theorem~\ref{linear-decay}. Set $\wt h_1 = h_1 - (h_1, \varphi_0)_{S_R} \varphi_0$. Notice that it holds $(\wt h_1, \varphi_0)_{S_R} = 0$. Let
\begin{equation*}
(\rho_{2 +}, \bu_{2 +}, \bu_{2 +}, h_2) (\cdot, s) = \int_0^s e^{- \wt A_q (s - r)} (- \lambda_1 \rho_{1 +}, - \lambda_1 \bu_{1 +}, - \lambda_1 \bu_{1 -}, - \lambda_1 \wt h_1) (\cdot, r) \,\mathrm{d} r,
\end{equation*}
and then by the Duhamel principle, we see that $(\rho_{2 +}, \bu_{2 +}, \bu_{2 +}, \wt h_2)$ satisfies 
\begin{equation}
\label{eq-shiftshift}
\left\{\begin{aligned}
\pd_t \rho_{2 +} + \rho_{* +} \dv \bu_{2 +} & = - \lambda_1 \rho_{1 +} &\enskip &\text{ in $\Omega_+ \times(0, T)$}, \\
\rho_{* -} \dv \bu_{2 -} & = 0 &\enskip &\text{ in $\Omega_- \times (0, T)$}, \\
\rho_{* +} \pd_t \bu_{2 +} - \DV \bT_+ (\rho_{2 +}, \bu_{2 +}) & = - \rho_{* +} \lambda_1 \bu_{1 +} &\enskip &\text{ in $\Omega_+ \times (0, T)$}, \\
\rho_{* -} \pd_t \bu_{2 -} - \DV \bT_- (\bu_{2 -}, \mathsf{P} (\bu_{2 -}, \wt h_2)) & = - \rho_{* +} \lambda_1 \bu_{1 -} &\enskip &\text{ in $\Omega_- \times (0, T)$}, \\
\pd_t \wt h_2 - \frac{1}{\rho_{* -} - \rho_{* +}} \lbrac \langle \rho_* \bu_2, \bn \rangle \rbrac + \CM \bu_2 & = - \lambda_1 \wt h_1 &\enskip &\text{ on $S_R \times (0, T)$}, \\
\bB (\rho_{2 +}, \bu_{2 +}, \bu_{2 -}, \mathsf{P} (\bu_{2 -}, \wt h_2), \wt h_2) & = 0 &\enskip &\text{ on $S_R \times (0, T)$}, \\
\bu_{2 +} = 0, \qquad \langle \nabla \rho_{2 +}, \bn_+ \rangle & = 0 &\enskip &\text{ on $\Gamma_+ \times (0, T)$}, \\
(\rho_{2 +}, \bu_2, \wt h_2) \rvert_{t = 0} & = (0, 0, 0) &\enskip &\text{ on $\Omega_+ \times \dot \Omega \times S_R$},			
\end{aligned}\right.		
\end{equation}
From Theorem~\ref{th-semigroup}, we have
\begin{align*}
& \lVert (\rho_{2 +}, \bu_{2 +}, \bu_{2 +}, \wt h_2) (\cdot, s) \rVert_{\mathsf{B}_q} \\
& = \int_0^s e^{- \eta_1 (s - r)} \lVert (\rho_{1 +}, \bu_{1 +}, \bu_{1 -}, \wt h_1) (\cdot, r) \rVert_{\mathsf{B}_q} \,\mathrm{d} r \\
& \le C \bigg(\int_0^s e^{- \eta_1 (s - r)} \,\mathrm{d} r\bigg)^{1 \slash p'} \bigg(\int_0^s e^{- \eta_1 (s - r)} \lVert (\rho_{1 +}, \bu_{1 +}, \bu_{1 -}, \wt h_1) (\cdot, r) \rVert_{\mathsf{B}_q}^p \bigg)^{1 \slash p}.
\end{align*}
Choosing $\eta > 0$ suitably small if necessary, we may suppose that $0 < \eta p < \eta_1$ without loss of generality. Hence, for any $t \in (0, T)$ it holds
\begin{align*}
& \int_0^t (e^{\eta s} \lVert (\rho_{2 +}, \bu_{2 +}, \bu_{2 +}, \wt h_2) (\cdot, s) \rVert_{\mathsf{B}_q})^p \,\mathrm{d} s \\
& \le C \int_0^t \bigg(\int_0^s e^{\eta s p} e^{- \eta_1 (s - r)} \lVert (\rho_{1 +}, \bu_{1 +}, \bu_{1 -}, h_1) (\cdot, r) \rVert_{\mathsf{B}_q}^p \mathrm{d} r \bigg) \ds \\
& = C \int_0^t \bigg(\int_0^s e^{- (\eta_1 - p \eta) (s - r)} \Big(e^{\eta r} \lVert (\rho_{1 +}, \bu_{1 +}, \bu_{1 -}, h_1) (\cdot, r) \rVert_{\mathsf{B}_q} \Big)^p \mathrm{d} r \bigg) \ds \\
& = C \int_0^t \Big(e^{\eta r} \lVert (\rho_{1 +}, \bu_{1 +}, \bu_{1 -}, h_1) (\cdot, r) \rVert_{\mathsf{B}_q} \Big)^p \bigg(\int_r^t e^{- (\eta_1 - p \eta) (s - r)} \ds \bigg) \,\mathrm{d} r \\
& \le C (\eta_1 - p \eta)^{- 1} \int_0^t \Big(e^{\eta r} \lVert (\rho_{1 +}, \bu_{1 +}, \bu_{1 -}, h_1) (\cdot, r) \rVert_{\mathsf{B}_q} \Big)^p \,\mathrm{d} r.
\end{align*}
Combined with \eqref{shift-decay}, we have
\begin{equation}
\label{shift-decay2}
\lVert e^{\eta s} (\rho_{2 +}, \bu_{2 +}, \bu_{2 +}, \wt h_2) (\cdot, s) \rVert_{L_p (0, t; \mathsf{B}_q)} \le C \BJ_{p, q, T} (\bu_0, h_0, F_M, F_d, \bF_d, \bF_+, \bF_-, D, \bG; \eta)
\end{equation}
for any $t \in (0, T)$. If $(p_{2 +}, \bw_{2 +}, \bw_{2 +}, \wt \eta_2)$ satisfies the shifted equations:
\begin{equation*}
\left\{\begin{aligned}
\pd_t p_{2 +} + \lambda_1 p_{2 +} + \rho_{* +} \dv \bw_{2 +} & = \lambda_1 (\rho_{2 +} - \rho_{1 +}) &\enskip &\text{ in $\Omega_+ \times(0, T)$}, \\
\rho_{* -} \dv \bw_{2 -} & = 0 &\enskip &\text{ in $\Omega_- \times (0, T)$}, \\
\rho_{* +} \pd_t \bw_{2 +} + \rho_{* +} \lambda_1 \bw_{2 +} - \DV \bT_+ (p_{2 +}, \bw_{2 +}) & = \rho_{* +} \lambda_1 (\bu_{2 +} - \bu_{2 +}) &\enskip &\text{ in $\Omega_+ \times (0, T)$}, \\
\rho_{* -} \pd_t \bw_{2 -} + \rho_{* -} \lambda_1 \bw_{2 -} - \DV \bT_- (\bw_{2 -}, \mathsf{P} (\bw_{2 -}, \wt \eta_2)) & = \rho_{* -} \lambda_1 (\bu_{2 -} - \bu_{1 -}) &\enskip &\text{ in $\Omega_- \times (0, T)$}, \\
\pd_t \wt \eta_2 + \lambda_1 \wt \eta_2 - \frac{1}{\rho_{* -} - \rho_{* +}} \lbrac \langle \rho_* \bw_2, \bn \rangle \rbrac + \CM \bw_2 & = \lambda_1 (\wt h_2 - \wt h_1) &\enskip &\text{ on $S_R \times (0, T)$}, \\
\bB (p_{2 +}, \bw_{2 +}, \bw_{2 -}, \mathsf{P} (\bw_{2 -}, \wt \eta_2), \eta_2) & = 0 &\enskip &\text{ on $S_R \times (0, T)$}, \\
\bw_{2 +} = 0, \qquad \langle \nabla p_{2 +}, \bn_+ \rangle & = 0 &\enskip &\text{ on $\Gamma_+ \times (0, T)$}, \\
(p_{2 +}, \bw_2, \wt \eta_2) \rvert_{t = 0} & = (0, 0, 0) &\enskip &\text{ on $\Omega_+ \times \dot \Omega \times S_R$}.					
\end{aligned}\right.	
\end{equation*}	
we have
\begin{equation*}
\CJ_{p, q, T} (p_{2 +}, \bw_{2 +}, \bw_{2 -}, \wt h_2; \eta) \le C \BJ_{p, q, T} (\bu_0, h_0, F_M, F_d, \bF_d, \bF_+, \bF_-, D, \bG; \eta)
\end{equation*}
as follows from \eqref{shift-decay} and \eqref{shift-decay2}. Recalling that $(\rho_{2 +}, \bu_{2 +}, \bu_{2 +}, \wt h_2)$ solves \eqref{eq-shiftshift}, we have $\rho_{2 +} = p_{2 +}$, $\bw_{2 \pm} = \bu_{2 \pm}$, and $\wt \eta_2 = \wt h_2$ for $t \in (0, T)$. Hence, we obtain the estimate
\begin{equation}
\label{decay-decay}
\CJ_{p, q, T} (\rho_{2 +}, \bu_{2 +}, \bu_{2 +}, \wt h_2; \eta) \le C \BJ_{p, q, T} (\bu_0, h_0, F_M, F_d, \bF_d, \bF_+, \bF_-, D, \bG; \eta)
\end{equation}
Now, we let
\begin{equation*}
h_2 = \wt h_2 - \lambda_1 \int_0^t (h_1 (\cdot, s), \varphi_0)_{S_R} \ds \, \varphi_0.
\end{equation*}
In this case, the pressure term $\mathsf{P} (\bu_2, h_2)$ is given by
\begin{equation*}
\mathsf{P} (\bu_{2 -}, h_2) = \mathsf{P} (\bu_{2 -}, \wt h_2) - \frac{(N - 1) \lambda_1}{R^2} \int_0^t (h_1 (\cdot, s), \varphi_0)_{S_R} \ds \, \varphi_0
\end{equation*}
In fact, it holds
\begin{equation*}
\CA_{S_R} h = \CA_{S_R} \wt h + \frac{(N - 1) \lambda_1}{R^2} \int_0^t (h_1 (\cdot, s), \varphi_0)_{S_R} \ds \, \varphi_0,
\end{equation*}
and thus we have
\begin{equation*}
\mathsf{P} (\bu_{2 -}, h_2) - \Big(\mathsf{P} (\bu_{2 -}, \wt h_2) + \CC\Big) = 0 \qquad \text{on $S_R$},
\end{equation*}
where we abbreviate
\begin{equation*}
\CC = - \frac{(N - 1) \lambda_1}{R^2} \int_0^t (h_1 (\cdot, s), \varphi_0)_{S_R} \ds \, \varphi_0.
\end{equation*}
Notice that we have
\begin{equation*}
\Big(\nabla \Big(\mathsf{P} (\bu_{2 -}, h_2) - (\mathsf{P} (\bu_{2 -}, \wt h_2) + \CC)\Big), \nabla \psi \Big)_{\Omega_-} = 0 
\end{equation*}
for any $\psi \in \wh H^1_{q', 0} (\Omega_-)$ since $\varphi_0$ is a constant, i.e., $\nabla \varphi_0 = 0$. Since $(\rho_{2 +}, \bu_{2 +}, \bu_{2 +}, \wt h_2)$ satisfies \eqref{eq-shiftshift}, we see that $(\rho_{2 +}, \bu_{2 +}, \bu_{2 +}, \mathsf{P} (\bu_{2 -}, h_2), h_2)$ satisfies \eqref{eq-homo}. Besides, by \eqref{decay-decay} we obtain the estimate
\begin{equation}
\label{shift-decay3}
\lVert e^{\eta t} \pd_t (\rho_{2 +}, \bu_{2 +}, \bu_{2 +}, h_2) \rVert_{L_p (0, t; \mathsf{B}_q)} \le C \BJ_{p, q, T} (\bu_0, h_0, F_M, F_d, \bF_d, \bF_+, \bF_-, D, \bG; \eta)
\end{equation}
for any $t \in (0, T)$. To estimate $\lVert e^{\eta t} (\rho_{2 +}, \bu_{2 +}, \bu_{2 +}, h_2) \rVert_{L_p (0, T; \mathsf{D}_q)}$, we use the following lemma.
\begin{lemm}
\label{lemm-steady} Let $1 < q < \infty$. Let $\rho_+ \in H^3_q (\Omega_+)$, $\bu_+ \in H^2_q (\Omega_+)^N$, $\bu_- \in H^2_q (\Omega_-)^N \cap J_q (\Omega_-)$, and $h \in W^{3 - 1 \slash q}_q (S_R)$ satisfy
\begin{equation}
\label{eq-steady}
\left\{\begin{aligned}
\rho_{* +} \dv \bu_+ & = F_M &\enskip &\text{ in $\Omega_+$}, \\
- \DV \bT_+ (\rho_+, \bu_+) & = \bF_+ &\enskip &\text{ in $\Omega_+$}, \\
- \DV \bT_- (\bu_-, \mathsf{P} (\bu_-, h)) & = \bF_- &\enskip &\text{ in $\Omega_-$}, \\
- \frac{1}{\rho_{* -} - \rho_{* +}} \lbrac \langle \rho_* \bu, \bn \rangle \rbrac + \CM \bu & = D &\enskip &\text{ on $S_R$}, \\
\bB (\rho_+, \bu_+, \bu_-, \mathsf{P} (\bu_-, h)) & = 0 &\enskip &\text{ on $S_R$}, \\
\bu_+ = 0, \qquad \langle \nabla \rho_+, \bn_+ \rangle & = 0 &\enskip &\text{ on $\Gamma_+$}
\end{aligned}\right.			
\end{equation}
with $(F_M, \bF_+, \bF_-, D) \in \mathsf{B}_q$. Then, it holds
\begin{equation}
\label{est-steady}
\lVert (\rho_+, \bu_+ \bu_-, h) \rVert_{\mathsf{D}_q} \le C \Big(\lVert (F_M, \bF_+, \bF_-, D) \rVert_{\mathsf{B}_q} + \lvert (h, \varphi_0)_{S_R} \rvert\Big)
\end{equation}
with some constant $C$.
\end{lemm}
We will give the proof of Lemma~\ref{lemm-steady} in the next section and we continue to show Theorem~\ref{linear-decay}. From \eqref{eq-homo}, we see that $(\rho_{2 +}, \bu_{2 +}, \bu_{2 -}, h_2)$ satisfies the elliptic equations:
\begin{equation*}
\left\{\begin{aligned}
\rho_{* +} \dv \bu_{2 +} & = - \lambda_1 \rho_{1 +} - \pd_t \rho_{2 +} &\enskip &\text{ in $\Omega_+$}, \\
\rho_{* -} \dv \bu_{2 -} & = 0 &\enskip &\text{ in $\Omega_-$}, \\
- \DV \bT_+ (\rho_{2 +}, \bu_{2 +}) & = - \rho_{* +} \lambda_1 \bu_{1 +} - \rho_{* +} \pd_t \bu_{2 +} &\enskip &\text{ in $\Omega_+$}, \\
- \DV \bT_- (\bu_{2 -}, \mathsf{P} (\bu_{2 -}, h_2)) & = - \rho_{* +} \lambda_1 \bu_{1 -} - \rho_{* -} \pd_t \bu_{2 -} &\enskip &\text{ in $\Omega_-$}, \\
- \frac{1}{\rho_{* -} - \rho_{* +}} \lbrac \langle \rho_* \bu_2, \bn \rangle \rbrac + \CM \bu_2 & = - \lambda_1 h_1 - \pd_t h_2 &\enskip &\text{ on $S_R$}, \\
\bB (\rho_{2 +}, \bu_{2 +}, \bu_{2 -}, \mathsf{P} (\bu_{2 -}, h_2), h_2) & = 0 &\enskip &\text{ on $S_R$}, \\
\bu_{2 +} = 0, \qquad \langle \nabla \rho_{2 +}, \bn_+ \rangle & = 0 &\enskip &\text{ on $\Gamma_+$},
\end{aligned}\right.		
\end{equation*}
According to Lemma~\ref{lemm-steady} and \eqref{shift-decay3}, we obtain
\begin{equation*}
\begin{split}
& \lVert e^{\eta t} (\rho_{2 +}, \bu_{2 +}, \bu_{2 +}, h_2) \rVert_{L_p (0, T; \mathsf{D}_q)} \\
& \quad \le C \bigg\{\BJ_{p, q, T} (\bu_0, h_0, F_M, F_d, \bF_d, \bF_+, \bF_-, D, \bG; \eta) + \bigg(\int_0^T (e^{\eta s} \lvert (h_2 (\cdot, s), \varphi_0)_{S_R} \rvert)^p \ds\bigg)^{1 \slash p} \bigg\}.
\end{split}
\end{equation*}
Finally, we find that $\rho_+ = \rho_{1 +} + \rho_{2 +}$, $\bu_\pm = \bu_{1 \pm} + \bu_{2 \pm}$, $\pi_- = \pi_{1 -} + \mathsf{P} (\bu_{2 -}, h_2)$, and $h = h_1 + h_2$ satisfy \eqref{eq-linear}. Especially, from the estimate
\begin{align*}
& \bigg(\int_0^T (e^{\eta s} \lvert (h_2 (\cdot, s), \varphi_0)_{S_R} \rvert)^p \ds\bigg)^{1 \slash p} \\
& \quad \le \bigg(\int_0^T (e^{\eta s} \lvert (h (\cdot, s), \varphi_0)_{S_R} \rvert)^p \ds\bigg)^{1 \slash p} + \bigg(\int_0^T (e^{\eta s} \lvert (h_1 (\cdot, s), \varphi_0)_{S_R} \rvert)^p \ds\bigg)^{1 \slash p} \\
& \quad \le \bigg(\int_0^T (e^{\eta s} \lvert (h (\cdot, s), \varphi_0)_{S_R} \rvert)^p \ds\bigg)^{1 \slash p} + C \BJ_{p, q, T} (\bu_0, h_0, F_M, F_d, \bF_d, \bF_+, \bF_-, D, \bG; \eta),
\end{align*}
which follows from \eqref{shift-decay}, we see that $(\rho_+, \bu_+, \bu_-, \pi_{2 -}, h_2)$ satisfies the required estimate. This completes the proof of Theorem~\ref{linear-decay}. 

\section{Exponential stability of continuous analytic semigroup}
\label{sec-decay}
\noindent
In this section, we shall prove Theorem~\ref{th-semigroup} and Theorem~\ref{lemm-steady} given in the previous section. To prove Theorem~\ref{th-semigroup}, we consider the resolvent problem associated to \eqref{Cauchy}:
\begin{equation}
\label{resolvent}
(\lambda \bI - A_q) \CU = F
\end{equation}
for $F = (F_M, \bF_+, \bF_-, D) \in \dot{\mathsf{B}}_q$ and $\CU = (\rho_+, \bu_+, \bu_-, h) \in \mathsf{D} (\wt A_q) = \mathsf{D} (A_q) \cap \dot{\mathsf{B}}_q$. According to \cite[Sec.~6.1]{W20}, there exists $\varepsilon_0^* \in (0, \pi \slash 2)$ and $\lambda_0 > 0$ such that for any $\varepsilon_0 \in (0, \varepsilon_0^*)$ and $\lambda \in \Sigma_{\varepsilon_0, \lambda_0}$, the resolvent estimate
\begin{equation}
\label{resolvent-est}
\lvert \lambda \rvert \lVert \CU \rVert_{\mathsf{B}_q} + \lVert \CU \rVert_{\mathsf{D}_q} \le C \lVert F \rVert_{\mathsf{B}_q}
\end{equation}
holds. To prove Theorem~\ref{th-semigroup}, we shall prove that the resolvent set of $\wt A_q$ contains $\BC_+ := \{z \in \BC \,\colon\, \mathrm{Re}\, z \ge 0\}$. 
\begin{theo}
Let $q \in (1, \infty)$. Then, for any $\lambda \in \BC_+$ and $F = (F_M, \bF_+, \bF_-, D) \in \dot{\mathsf{B}}_q$, the problem \eqref{resolvent} admits a unique solution $\CU \in \mathsf{D} (\wt A_q)$ that satisfies the estimate \eqref{resolvent-est}.
\end{theo}
Recalling Lemma~\ref{lem-invariant}, we obtain the following lemma.
\begin{lemm}
\label{lemlem-resolvent}
Let $q \in (1, \infty)$. Then there exists $\varepsilon_0^* \in (0, \pi \slash 2)$ such that for any $\varepsilon_0 \in (0, \varepsilon_0^*)$ there is a constant $\lambda_0 > 0$ with the following property holds. For every $\lambda \in \Sigma_{\varepsilon_0, \lambda_0}$ and $F = (F_M, \bF_+, \bF_-, D) \in \dot{\mathsf{B}}_q$, the resolvent problem~\eqref{resolvent} has a unique solution $\CU = (\rho_+, \bu_+, \bu_-, h) \in \mathsf{D} (\wt A_q)$ possessing the estimate \eqref{resolvent-est}.
\end{lemm}
In view of Lemma~\ref{lemlem-resolvent}, it suffices to show the next theorem.
\begin{theo}
\label{semi-simple}
Let $q \in (1, \infty)$ and $\lambda_0$ be the same number as in Lemma~\ref{lemlem-resolvent}. For any $\lambda \in Q_{\lambda_0} := \{\lambda \in \BC \,\colon\, \mathrm{Re}\,\lambda \ge 0, \enskip \lvert \lambda \rvert \le \lambda_0\}$ and $F = (F_M, \bF_+, \bF_-, D) \in \dot{\mathsf{B}}_q$, the equation \eqref{resolvent} admits a unique solution $\CU = (\rho_+, \bu_+, \bu_-, h) \in \mathsf{D} (\wt A_q)$ that enjoys the estimate $\lVert (\rho_+, \bu_+, \bu_-, h) \rVert_{\mathsf{D}_q} \le C \lVert (F_M, \bF_+, \bF_-, D) \rVert_{\mathsf{B}_q}$ with some constant $C$ independent of $\lambda$.
\end{theo}
In the following, we shall prove Lemma~\ref{semi-simple}. We first observe that 
\begin{equation}
\label{embd}
A_q \mathsf{D} (\wt A_q) \subset \dot{\mathsf{B}}_q.
\end{equation}
In fact, for $(\rho_+, \bu_+, \bu_-, h) \in \mathsf{D} (\wt A_q)$ we set $A_q (\rho_+, \bu_+, \bu_-, h) = (F_M, \bF_+, \bF_-, D)$, i.e.,
\begin{equation}
\label{zero-resolvent}
\left\{\begin{aligned}
\rho_{* +} \dv \bu_+ & = F_M &\enskip &\text{ in $\Omega_+$}, \\
- \DV \bT_+ (\rho_+, \bu_+) & = \bF_+ &\enskip &\text{ in $\Omega_+$}, \\
- \DV \bT_- (\bu_-, \mathsf{P} (\bu_-, h)) & = \bF_- &\enskip &\text{ in $\Omega_-$}, \\
- \frac{1}{\rho_{* -} - \rho_{* +}} \lbrac \langle \rho_* \bu, \bn \rangle \rbrac + \CM \bu & = D &\enskip &\text{ on $S_R$}, \\
\bB (\rho_+, \bu_+, \bu_-, \mathsf{P} (\bu_-, h)) & = 0 &\enskip &\text{ on $S_R$}, \\
\bu_+ = 0, \qquad \langle \nabla \rho_+, \bn_+ \rangle & = 0 &\enskip &\text{ on $\Gamma_+$}.				
\end{aligned}\right.			
\end{equation}
For any $\varphi \in \wh H^1_{q', 0} (\Omega_-)$, by \eqref{pressure}, we have
\begin{equation*}
(\bF_-, \nabla \varphi)_{\Omega_-} = (- \DV \bT_- (\bu_-, \mathsf{P} (\bu_-, h)), \nabla \varphi)_{\Omega_-} = (- \rho_{* -} \nabla \dv \bu_-, \nabla \varphi)_{\Omega_-} = 0,
\end{equation*}
which implies $\bF_- \in J_q (\Omega_-)$. From \eqref{zero-resolvent}$_4$ and the divergence theorem, it holds $(D, 1)_{S_R} = 0$ due to $\dv \bu_- = 0$ in $\Omega_-$. Hence, we obtain \eqref{embd}. \par
In view of Lemma~\ref{lemlem-resolvent}, the inverse $(\lambda_0 \bI - A_q)^{- 1}$ exists as a bounded linear operator from $\dot{\mathsf{B}}_q$ onto $\mathsf{D} (\wt A_q)$. Then, the equation \eqref{resolvent} is rewritten as
\begin{align*}
(F_M, \bF_+, \bF_-, D) & = (\lambda \bI - A_q) (\rho_+, \bu_+, \bu_-, h) \\
& = (\lambda - \lambda_0) (\rho_+, \bu_+, \bu_-, h) + (\lambda_0 \bI - A_q) (\rho_+, \bu_+, \bu_-, h) \\
& = (\bI + (\lambda - \lambda_0) (\lambda_0 \bI - A_q)^{- 1}) (\lambda_0 \bI - A_q) (\rho_+, \bu_+, \bu_-, h).
\end{align*}
If the inverse of $(\bI + (\lambda - \lambda_0) (\lambda_0 \bI - A_q)^{- 1})$ exists as a bounded linear operator from $\dot{\mathsf{B}}_q$ onto itself, it holds
\begin{equation*}
(\rho_+, \bu_+, \bu_-, h) = (\lambda_0 \bI - A_q)^{- 1} \Big((\bI + (\lambda - \lambda_0) (\lambda_0 \bI - A_q)^{- 1})\Big)^{- 1} (F_M, \bF_+, \bF_-, D).
\end{equation*}
Hence, it remains to prove the existence of the operator $(\bI + (\lambda - \lambda_0) (\lambda_0 \bI - A_q)^{- 1})^{- 1}$. Here, $(\lambda_0 \bI - A_q)^{- 1}$ is a compact operator from $\dot{\mathsf{B}}_q$ onto itself due to the Rellich compact embedding theorem. Hence, in view of the Riesz-Schauder theory, it suffices to show that the kernel of the map $\bI + (\lambda - \lambda_0) (\lambda_0 \bI - A_q)^{- 1}$ is trivial, i.e., if $(F_M, \bF_+, \bF_-, D) \in \dot{\mathsf{B}}_q$ satisfies
\begin{equation}
\label{trivial}
\Big((\bI + (\lambda - \lambda_0) (\lambda_0 \bI - A_q)^{- 1})\Big) (F_M, \bF_+, \bF_-, D) = (0, 0, 0, 0),
\end{equation}
then $(F_M, \bF_+, \bF_-, D) = (0, 0, 0, 0)$. From \eqref{trivial}, it holds 
\begin{equation*}
(F_M, \bF_+, \bF_-, D) = - (\lambda - \lambda_0) (\lambda_0 \bI - A_q)^{- 1} (F_M, \bF_+, \bF_-, D) \in \mathsf{D} (\wt A_q)
\end{equation*}
so that we observe $(\lambda_0 \bI - A_q) (F_M, \bF_+, \bF_-, D) = - (\lambda - \lambda_0) (F_M, \bF_+, \bF_-, D)$, i.e., $(\lambda \bI - A_q) (F_M, \bF_+, \bF_-, D) = (0, 0, 0, 0)$. This equivalents to the fact that $(\rho_+, \bu_+, \bu_-, h) \in \mathsf{D} (\wt A_q)$ satisfies the homogeneous equations:
\begin{equation}
\label{eq-resolvent}	
\left\{\begin{aligned}
\lambda \rho_+ + \rho_{* +} \dv \bu_+ & = 0 &\enskip &\text{ in $\Omega_+$}, \\
\rho_{* -} \dv \bu_- & = 0 &\enskip &\text{ in $\Omega_-$}, \\
\rho_{* +} \lambda \bu_+ - \DV \bT_+ (\rho_+, \bu_+) & = 0 &\enskip &\text{ in $\Omega_+$}, \\
\rho_{* -} \lambda \bu_- - \DV \bT_- (\bu_-, \mathsf{P} (\bu_-, h)) & = 0 &\enskip &\text{ in $\Omega_-$}, \\
\lambda h - \frac{1}{\rho_{* -} - \rho_{* +}} \lbrac \langle \rho_* \bu, \bn \rangle \rbrac + \CM \bu & = 0 &\enskip &\text{ on $S_R$}, \\
\bB (\rho_+, \bu_+, \bu_-, \mathsf{P} (\bu_-, h)) & = 0 &\enskip &\text{ on $S_R$}, \\
\bu_+ = 0, \qquad \langle \nabla \rho_+, \bn_+ \rangle & = 0 &\enskip &\text{ on $\Gamma_+$}.
\end{aligned}\right.		
\end{equation}
We first notice that the spectrum of $\wt A_q$ is independent of $q$, so that we may let $q = 2$, cf., \cite{A94}. Taking the inner product of the problem for $\rho_+$ with $\rho_{* +}^{- 1} \gamma_{* +} \rho_+$ and $\rho_{* +} \kappa_+ \nabla \rho_+$, an integration by parts infers
\begin{equation}
\label{inner-product}
\begin{split}
\frac{\gamma_{* +}}{\rho_{* +}} \lambda \lVert \rho_+ \rVert_{L_2 (\Omega_+)}^2 + (\gamma_{* +} \rho_+, \dv \bu_+)_{\Omega_+} & = 0, \\
\kappa_+ \lambda \lVert \nabla \rho_+ \rVert_{L_2 (\Omega_+)}^2 + (\rho_{* +} \kappa_+ \nabla \rho_+, \nabla \dv \bu_+)_{\Omega_+} & = 0.
\end{split}
\end{equation}
On the other hand, the inner product of the equations for $\bu_\pm$ with $\bu_\pm$ by an integration by parts leads to
\begin{equation}
\label{innerproducts}
\begin{split}
0 & = \sum_{\ell = \pm} \bigg(\lambda \rho_{* \ell} \lVert \bu_\ell \rVert_{L_2 (\Omega_\ell)}^2 + \frac{\mu_\ell}{2} \lVert \bD (\bu_\ell) \rVert_{L_2 (\Omega_\ell)}^2 \bigg) - (\bT_- (\bu_-, \mathsf{P} (\bu_-, h)) \bn \vert_- , \bu_-)_{S_R} \\
& \quad + (\bT_+ (\rho_+, \bu_+) \bn \vert_+ , \bu_+)_{S_R} + (\nu_+ - \mu_+) \lVert \dv \bu_+ \rVert_{L_2 (\Omega_+)}^2 \\
& \quad - (\gamma_{* +} \rho_+, \dv \bu_+)_{\Omega_+} - (\rho_{* +} \kappa_+ \nabla \rho_+, \dv \bu_+)_{\Omega_+}.
\end{split}
\end{equation}
Since we have
\begin{align*}
(\bT_- (\bu_-, \mathsf{P} (\bu_-, h)) \bn \vert_- , \bu_-)_{S_R} & = (\bT_- (\bu_-, \mathsf{P} (\bu_-, h)) \bn \vert_- , \Pi_\bn \bu_- + \langle \bu_-, \bn \rangle \bn)_{S_R}, \\
(\bT_+ (\rho_+, \bu_+) \bn \vert_+ , \bu_+)_{S_R} & = (\bT_+ (\rho_+, \bu_+) \bn \vert_+ , \Pi_\bn \bu_+ + \langle \bu_+, \bn \rangle \bn)_{S_R},
\end{align*}
it holds
\begin{align*}
& - (\bT_- (\bu_-, \mathsf{P} (\bu_-, h)) \bn \vert_- , \bu_-)_{S_R} + (\bT_+ (\rho_+, \bu_+) \bn \vert_+ , \bu_+)_{S_R} \\
& \quad = - (\langle \bT_- (\bu_-, \mathsf{P} (\bu_-, h)) \bn \vert_-, \bn \rangle \bn, \langle \bu_-, \bn \rangle \bn)_{S_R} + (\langle \bT_+ (\rho_+, \bu_+) \bn \vert_+, \bn \rangle \bn, \langle \bu_+, \bn \rangle \bn)_{S_R}.
\end{align*}
Moreover, by
\begin{align*}
(\langle \bT_- (\bu_-, \mathsf{P} (\bu_-, h)) \bn \vert_-, \bn \rangle \bn, \langle \bu_-, \bn \rangle \bn)_{S_R} & = \bigg(- \frac{\rho_{* -} \sigma}{\rho_{* -} - \rho_{* +}} \CA_{S_R} h + \sum_{j = 1}^N (h, \varphi_j)_{S_R} \varphi_j, \langle \bu_-, \bn \rangle \bn \bigg)_{S_R}, \\
(\langle \bT_+ (\rho_+, \bu_+) \bn \vert_+, \bn \rangle \bn, \langle \bu_+, \bn \rangle \bn)_{S_R} & = \bigg(- \frac{\rho_{* +} \sigma}{\rho_{* -} - \rho_{* +}} \CA_{S_R} h, \langle \bu_+, \bn \rangle \bn \bigg)_{S_R},
\end{align*}
the equation of $h$ implies
\begin{align*}
& - (\langle \bT_- (\bu_-, \mathsf{P} (\bu_-, h)) \bn \vert_-, \bn \rangle \bn, \langle \bu_-, \bn \rangle \bn)_{S_R} + (\langle \bT_+ (\rho_+, \bu_+) \bn \vert_+, \bn \rangle \bn, \langle \bu_+, \bn \rangle \bn)_{S_R} \\
& \quad = \bigg(\sigma \CA_{S_R} h + \frac{\rho_{* +} - \rho_{* -}}{\rho_{* -} \sigma} \sum_{j = 1}^N (h, \varphi_j)_{S_R} \varphi_j,  \frac{\lbrac \langle \rho_* \bu, \bn \rangle \rbrac}{\rho_{* -} - \rho_{* +}} \bigg)_{S_R} \\
& \quad = \bigg(\sigma \CA_{S_R} h + \frac{\rho_{* +} - \rho_{* -}}{\rho_{* -} \sigma} \sum_{j = 1}^N (h, \varphi_j)_{S_R} \varphi_j, \lambda h + \CM \bu \bigg)_{S_R}.
\end{align*}
Since the components of $\omega = (\omega_1, \dots, \omega_N)$ are eigenfunctions of the Laplace-Beltrami operator $\Delta_{S_R}$, it holds $\Delta_{S_R} \omega_j = 0$, $j = 1, \dots, N$, so that
\begin{equation}
\label{laplace}
\begin{split}
& - (\bT_- (\bu_-, \mathsf{P} (\bu_-, h)) \bn \vert_- , \bu_-)_{S_R} + (\bT_+ (\rho_+, \bu_+) \bn \vert_+ , \bu_+)_{S_R} \\
& = \ov \lambda (\sigma \CA_{S_R} h, h)_{S_R} + \frac{\rho_{* +} - \rho_{* -}}{\rho_{* -} \sigma} \sum_{j = 1}^N \ov \lambda (h, \varphi_j)_{S_R}^2.
\end{split}
\end{equation}
From \eqref{inner-product} and \eqref{laplace}, we see that \eqref{innerproducts} can be rewritten as
\begin{equation}
\label{lambdaproduct}
\begin{split}
0 & = \lambda \bigg(\rho_{* +} \lVert \bu_+ \rVert_{L_2 (\Omega_+)}^2 + \rho_{* -} \lVert \bu_- \rVert_{L_2 (\Omega_-)}^2 + \frac{\gamma_{* +}}{\rho_{* +}} \lVert \rho_+ \rVert_{L_2 (\Omega_+)}^2 + \kappa_+ \lVert \nabla \rho_+ \rVert_{L_2 (\Omega_+)}^2 \bigg) \\
& \quad + \ov \lambda (\sigma \CA_{S_R} h, h)_{S_R} + \frac{\rho_{* +} - \rho_{* -}}{\rho_{* -} \sigma} \sum_{j = 1}^N \ov \lambda (h, \varphi_j)_{S_R}^2 + \frac{\mu_+}{2} \lVert \bD (\bu_+) \rVert_{L_2 (\Omega_+)}^2 + \frac{\mu_-}{2} \lVert \bD (\bu_-) \rVert_{L_2 (\Omega_-)}^2 \\
& \quad  + (\nu_+ - \mu_+) \lVert \dv \bu_+ \rVert_{L_2 (\Omega_+)}^2.
\end{split}
\end{equation}
Hence, taking the real part of \eqref{lambdaproduct}, it follows that
\begin{equation}
\label{lam}
\begin{split}
0 & = (\mathrm{Re}\, \lambda) \bigg(\rho_{* +} \lVert \bu_+ \rVert_{L_2 (\Omega_+)}^2 + \rho_{* -} \lVert \bu_- \rVert_{L_2 (\Omega_-)}^2 + \frac{\gamma_{* +}}{\rho_{* +}} \lVert \rho_+ \rVert_{L_2 (\Omega_+)}^2 + \kappa_+ \lVert \nabla \rho_+ \rVert_{L_2 (\Omega_+)}^2 \\
& \qquad\qquad + (\sigma \CA_{S_R} h, h)_{S_R} + \frac{\rho_{* +} - \rho_{* -}}{\rho_{* -} \sigma} \sum_{j = 1}^N (h, \varphi_j)_{S_R}^2 \bigg) \\
& \quad + \frac{\mu_+}{2} \lVert \bD (\bu_+) \rVert_{L_2 (\Omega_+)}^2 + \frac{\mu_-}{2} \lVert \bD (\bu_-) \rVert_{L_2 (\Omega_-)}^2 + (\nu_+ - \mu_+) \lVert \dv \bu_+ \rVert_{L_2 (\Omega_+)}^2.
\end{split}		
\end{equation}
To handle $(\sigma \CA_{S_R} h, h)_{S_R}$, we introduce the following lemma essentially proved in \cite[Lem.~4.5]{S18}.
\begin{lemm}
\label{laplacebeltrami}
Let $\CA_{S_R}$ be defined on $L^2 (S_R)$ with domain $H^{2, 2} (S_R)$. Then the following holds.
\begin{enumerate}
\item $\CA_{S_R}$ is self-adjoint. The spectrum of $\CA_{S_R}$ consists entirely of eigenvalues of finite algebraic multiplicity and is given by $\sigma (\CA_{S_R}) = \{R^{- 2} (k (k + 1) - 2) \,\colon \, k \in \BN_0\}$.
\item There is precisely one negative eigenvalue $- 2 \slash R^2$ with eigenfunction $e \equiv 1$, which is simple.
\item The kernel of $\CA_{S_R}$ is spanned by $\varphi_j$, $j = 1, \dots, N$.
\item $\CA_{S_R}$ is positive semi-definite on $L^2_0 (S_R) := \{\eta \in L^2 (S_R) \, \colon \, (\eta, 1)_{S_R} = 0\}$ and positive definite on
\begin{equation*}
L^2_0 (S_R) \cap \mathsf{R} (\CA_{S_R}) = \{\eta \in L^2 (S_R) \,\colon\, (\eta, \varphi_j)_{S_R} = 0, \enskip j = 0, \dots, N \}.
\end{equation*}
\end{enumerate}	
\end{lemm}
From the equation of $h$ and the divergence theorem of Gauss, we have
\begin{align*}
0 & = \bigg(\lambda h - \frac{1}{\rho_{* -} - \rho_{* +}} \lbrac \langle \rho_* \bu, \bn \rangle \rbrac + \CM \bu, 1 \bigg)_{S_R} \\
& = \lambda \int_{S_R} h \,\mathrm{d} \omega - \frac{\rho_{* -}}{\rho_{* -} - \rho_{* +}} (\dv \bu_-, 1)_{\Omega_-},
\end{align*}
where $\lvert S_R \rvert$ denotes the area of $S_R$. Since $\dv \bu_- = 0$ in $\Omega_-$, we have
\begin{equation}
\label{cond-h1}
\int_{S_R} h \,\mathrm{d} \omega = 0
\end{equation}
provided that $\mathrm{Re}\, \lambda \ge 0$ with $\lambda \ne 0$. Hence, it follows from Lemma~\ref{laplacebeltrami} that $(\sigma \CA_{S_R} h, h)_{S_R}$ is positive semi-definite. Besides, noting that $\lVert \bD (\bu_+) \rVert_{L_2 (\Omega_+)}^2 \ge (4 \slash N) \lVert \dv \bu_+\rVert_{L_2 (\Omega_+)}^2$, we have
\begin{equation}
\label{lambdalambda}
\begin{split}
0 & \ge (\mathrm{Re}\, \lambda) \bigg(\rho_{* +} \lVert \bu_+ \rVert_{L_2 (\Omega_+)}^2 + \rho_{* -} \lVert \bu_- \rVert_{L_2 (\Omega_-)}^2 + \frac{\gamma_{* +}}{\rho_{* +}} \lVert \rho_+ \rVert_{L_2 (\Omega_+)}^2 + \kappa_+ \lVert \nabla \rho_+ \rVert_{L_2 (\Omega_+)}^2 \bigg) \\
& \quad + \frac{\mu_-}{2} \lVert \bD (\bu_-) \rVert_{L_2 (\Omega_-)}^2 + \bigg(\nu_+ - \frac{N - 2}{N} \mu_+ \bigg) \lVert \dv \bu_+ \rVert_{L_2 (\Omega_+)}^2.
\end{split}
\end{equation}
for all $\lambda$ such that $\mathrm{Re}\, \lambda \ge 0$, $\lambda \ne 0$. Since $\nu_+ \ge ((N - 2) \slash N) \mu_+ \ge 0$, we have $\rho_+ = \bu_+ = \bu_- = 0$ when $\mathrm{Re}\,\lambda > 0$ and $\lambda \ne 0$. Hence, from the equation of $h$, we obtain $h = 0$ if $\lambda > 0$. If $\mathrm{Re}\, \lambda = 0$, the inequality \eqref{lambdalambda} yields $\bD (\bu_-) = 0$ in $\Omega_-$ and $\dv \bu_+ = 0$ in $\Omega_+$. Recalling~\eqref{lam}, we find that $\bD (\bu_+) = 0$ in $\Omega_+$. Using the Korn inequality, we observe $\bu_+ = \bu_- = 0$ due to the no-slip boundary condition on $\Gamma_+$ and the boundary condition $\Pi_\bn \bu_- \vert_- - \Pi_\bn \bu_+ \vert_+ = 0$ on $S_R$, see also \cite[Lem.~1.2.1]{PS16}. Then, from the equation of $\rho_+$, we find that $\rho_+ = 0$ because $\lambda \ne 0$. Besides, by the equation of $h$, we also obtain $h = 0$ due to $\lambda \ne 0$. This investigation shows that $\lambda$ is not an eigenvalue of $\wt A_q$ if $\mathrm{Re}\,\lambda \ge 0$ and $\lambda \ne 0$. \par
We now show that $\lambda = 0$ belongs to a resolvent set of $\wt A_q$ as well. As we discussed above, we easily observe that $\bu_+ = \bu_- = 0$. By \eqref{eq-resolvent}$_4$, we see that $\mathsf{P} (\bu_-, h)$ is a constant in $\Omega_-$. Here, by the interface condition for the stress tensor, we have 
\begin{equation*}
\mathsf{P} (\bu_-, h) \vert_- + \frac{\rho_{* -} \sigma}{\rho_{* -} + \rho_{* +}} \CA_{S_R} h + \sum_{j = 1}^N (h, \varphi_j)_{S_R} \varphi_j = 0 \qquad \text{on $S_R$}.
\end{equation*}
Integrating this formula on $S_R$ and using~\eqref{cond-h1}, we arrive at $\mathsf{P} (\bu_-, h) = 0$ on $S_R$, i.e.,
\begin{equation*}
\frac{\rho_{* -} \sigma}{\rho_{* -} + \rho_{* +}} \CA_{S_R} h + \sum_{j = 1}^N (h, \varphi_j)_{S_R} \varphi_j = 0 \qquad \text{on $S_R$}.
\end{equation*}
Taking the inner product of this identity with $\varphi_j$, we observe $(h, \varphi_j)_{S_R}= 0$ due to Lemma~\ref{laplacebeltrami}. Hence, we see that $\CA_{S_R} h = 0$ on $S_R$, which implies $h = 0$ on $S_R$. Now, from \eqref{eq-resolvent}$_3$, it holds $\gamma_{* +} \nabla \rho_+ - \rho_{* +} \kappa_+ \nabla \Delta \rho_+ = 0$ in $\Omega_+$, where we have $\langle \nabla \rho_+, \bn \rangle = 0$ on $S_R$. Taking the inner product of this elliptic problem with $\nabla \rho_+$, we have
\begin{equation*}
\gamma_{* +} \lVert \nabla \rho_+ \rVert_{L_2 (\Omega_+)}^2 + \rho_{* +} \kappa_+ \lVert \Delta \rho_+ \rVert_{L_2 (\Omega_+)}^2 = 0.
\end{equation*}
This gives that $\rho_+$ is a constant. However, recalling $\CA_{S_R} h = 0$ on $S_R$ and the stress boundary condition, we obtain $- \gamma_{* +} \rho_+ = 0$ in $\Omega_+$. Since $\gamma_{* +} > 0$, we deduce that $\rho_+ = 0$ in $\Omega_+$. This completes the proof of Theorem~\ref{semi-simple}. \par
Finally, we give the proof of Theorem~\ref{lemm-steady}. Let $\wt h = h - (h, \varphi_0)_{S_R} \varphi_0$. Then, it holds
\begin{equation}
(\nabla \mathsf{P} (\bu_-, \wt h), \nabla \psi)_{\Omega_-} = (\DV (\mu_- \bD (\bu_-)) - \rho_{* -} \nabla \dv \bu_-, \nabla \psi)_{\Omega_-} = (\nabla \mathsf{P} (\bu_-, h), \nabla \psi)_{\Omega_-}
\end{equation}
for any $\psi \in \wh H^1_{q', 0} (\Omega_-)$ subject to
\begin{align*}
\mathsf{P} (\bu_-, \wt h) & = \mu_- \langle \bD (\bu_-) \bn, \bn \rangle + \frac{\rho_{* -} \sigma}{\rho_{* -} - \rho_{* +}} \CA_{S_R} h + \frac{\sigma (N - 1)}{R^2} (h, \varphi_0)_{S_R} \varphi_0 \\
& \quad + \sum_{j = 1}^N (h, \varphi_j)_{S_R} \varphi_j - \rho_{* -} \dv \bu_- \qquad\qquad\qquad\qquad\qquad\qquad\qquad\qquad \text{on $S_R$}.
\end{align*}
Namely, we have
\begin{equation*}
\mathsf{P} (\bu_-, \wt h) = \mathsf{P} (\bu_-, h) + \frac{\sigma (N - 1)}{R^2} (h, \varphi_0)_{S_R} \varphi_0.
\end{equation*}
Noting that $\varphi_0$ is a constant, i.e., $\nabla \varphi_0 = 0$, we observe that $(\rho_+, \bu_+, \bu_-, \wt h)$ satisfies \eqref{eq-steady} with $(\rho_+, \bu_+, \bu_-, \wt h) \in \mathsf{D} (\wt A_q)$. Hence, by \eqref{embd} and Theorem~\ref{semi-simple} with $\lambda = 0$, it holds $\lVert (\rho_+, \bu_+, \bu_-, \wt h) \rVert_{\mathsf{D}_q} \le C \lVert (F_M, \bF_+, \bF_-, D) \rVert_{\mathsf{B}_q}$. Therefore, combined with the estimate:
\begin{equation*}
\lVert (\rho_+, \bu_+ \bu_-, h) \rVert_{\mathsf{D}_q} \le C \Big(\lVert (\rho_+, \bu_+ \bu_-, \wt h) \rVert_{\mathsf{D}_q} + \lvert (h, \varphi_0)_{S_R} \rvert\Big),
\end{equation*}
we obtain \eqref{est-steady}. This completes the proof of Lemma~\ref{lemm-steady}.

\section{Nonlinear well-posedness}
\label{sec-nonlinear}
\subsection{Local well-posedness}
Before we prove Theorem~\ref{th-main}, we state the local well-posedness result of \eqref{eq-1.5}.
\begin{theo}
\label{th-local}
Let $2 < p < \infty$, $N < q < \infty$, $2 \slash p + N \slash q < 1$. Besides, let $T > 0$. Suppose that Assumptions~\ref{asp-1} and \ref{asp-2} hold. Then, there exists a number $\varepsilon > 0$ depending on $T$ such that if initial data $\rho_{0 +} \in B^{3 - 2 \slash p}_{q, p} (\Omega_+)$, $\bu_{0 \pm} \in B^{2 (1 - 1 \slash p)}_{q, p} (\Omega_\pm)$, and $h_0 \in B^{3 - 1 \slash p - 1 \slash q}_{q, p} (S_R)$ satisfies the smallness condition:
\begin{equation*}
\lVert \rho_{0 +} \rVert_{B^{3 - 2 \slash p}_{q, p} (\Omega_+)} + \lVert \bu_{0 +} \rVert_{B^{2 (1 - 1 \slash p)}_{q, p} (\Omega_+)} + \lVert \bu_{0 -} \rVert_{B^{2 (1 - 1 \slash p)}_{q, p} (\Omega_-)} + \lVert h_0 \rVert_{B^{3 - 1 \slash p - 1 \slash q}_{q, p} (S_R)} \le \varepsilon
\end{equation*}
and the compatibility conditions:
\begin{equation*}
\left\{\begin{aligned}
\rho_{* -} \dv \bu_{0 -} = f_d (\bu_{0 -}, h_0) & = \rho_{* -} \dv \bff_d (\bu_{0 -}, h_0) & \quad & \text{in $\Omega_-$}, \\
\Pi_\bn (\mu_- \bD (\bu_-) \bn) \rvert_- - \Pi_\bn (\mu_+ \bD(\bu_+) \bn) \rvert_+ & = g (\rho_{0 +}, \bu_{0 +}, \bu_{0 -}, h_0) & \quad & \text{on $S_R$}, \\
\Pi_\bn \bu_- \rvert_- - \Pi_\bn \bu_+ \rvert_+ & = \bh (\bu_{0 +}, \bu_{0 -}, h_0) & \quad & \text{on $S_R$}, \\
\langle \nabla \rho_{0 +}, \bn \rangle \rvert_+ & = k_- (\rho_{0 +}, h_0) & \quad & \text{on $S_R$},	\\
\langle \nabla \rho_{0 +}, \bn_+ \rangle = 0, \quad \bu_{0 +} & = 0	& \quad & \text{on $\Gamma_+$},
\end{aligned}\right.
\end{equation*}
the problem \eqref{eq-1.5} has a unique solution $(\rho_+, \bu_+, \bu_-, \pi_-, h) \in \CS_{p, q} (0, T)$ and the estimate
\begin{align*}
& \CJ_{p, q, T} (\rho_+, \bu_+, \bu_-, \pi_-, h; 0) + \lVert \rho_+ \rVert_{L_\infty (0, T; B^{3 - 2 \slash p}_{q, p} (\Omega_+))} \\
& \qquad + \sum_{\ell = \pm} \lVert \bu_\ell \rVert_{L_\infty (0, T; B^{2 (1 - 1 \slash p)}_{q, p} (\Omega_\ell))} + \lVert h \rVert_{L_\infty (0, T; B^{3 - 1 \slash p - 1 \slash q}_{q, p} (S_R))} \le \varepsilon.
\end{align*}
\end{theo}

Moving the lower-order terms $\sum_{j = 1}^N (h, \varphi_j)_{S_R} \varphi_j$ and $\CM \bu$ to the right-hand side and employing the similar argument to that in the proof of \cite[Thm.~3.7]{W20}, we can prove Theorem~\ref{th-local}, and so we may omit the details. Here, by $2 \slash p + N \slash q < 1$ it holds $B^{3 - 1 \slash p - 1 \slash q}_{q, p} (S_R) \hookrightarrow H^2_\infty (S_R)$, so that the condition \eqref{cond-psi} holds if $\varepsilon$ is so small.

\subsection{Global well-posedness}
Finally, we prove Theorem~\ref{th-main}. Assume that the initial data $\rho_{0 +} \in B^{3 - 2 \slash p}_{q, p} (\Omega_+)$, $\bu_{0 \pm} \in B^{2 (1 - 1 \slash p)}_{q, p} (\Omega_\pm)$, and $h_0 \in B^{3 - 1 \slash p - 1 \slash q}_{q, p} (S_R)$ satisfy the smallness condition
\begin{equation*}
\lVert \rho_{0 +} \rVert_{B^{3 - 2 \slash p}_{q, p} (\Omega_+)} + \lVert \bu_{0 +} \rVert_{B^{2 (1 - 1 \slash p)}_{q, p} (\Omega_+)} + \lVert \bu_{0 -} \rVert_{B^{2 (1 - 1 \slash p)}_{q, p} (\Omega_-)} + \lVert h_0 \rVert_{B^{3 - 1 \slash p - 1 \slash q}_{q, p} (S_R)} \le \varepsilon
\end{equation*}
with small constant $\varepsilon > 0$ as well as the compatibility conditions \eqref{compati-h0} and \eqref{compat-1}. In the following, we write
\begin{align*}
& \CJ := \lVert \rho_{0 +} \rVert_{B^{3 - 2 \slash p}_{q, p} (\Omega_+)} + \lVert \bu_{0 +} \rVert_{B^{2 (1 - 1 \slash p)}_{q, p} (\Omega_+)} + \lVert \bu_{0 -} \rVert_{B^{2 (1 - 1 \slash p)}_{q, p} (\Omega_-)} + \lVert h_0 \rVert_{B^{3 - 1 \slash p - 1 \slash q}_{q, p} (S_R)}, \\
& E_{p, q, T} (\rho_+, \bu_+, \bu_-, \pi_-, h; \eta) := \CJ_{p, q, T} (\rho_+, \bu_+, \bu_-, \pi_-, h; \eta) + \lVert \rho_+ \rVert_{L_\infty (0, T; B^{3 - 2 \slash p}_{q, p} (\Omega_+))} \\
& \qquad \qquad \qquad \qquad \qquad \qquad \quad + \sum_{\ell = \pm} \lVert \bu_\ell \rVert_{L_\infty (0, T; B^{2 (1 - 1 \slash p)}_{q, p} (\Omega_\ell))} + \lVert h \rVert_{L_\infty (0, T; B^{3 - 1 \slash p - 1 \slash q}_{q, p} (S_R))}
\end{align*}
for short. From the proof of \cite[Lem.~5.4]{S18} (cf. \cite[(3.212), (3.213)]{SS20}), there exists a constant $C$ independent of $T$ such that the estimate
\begin{equation}
\label{402}
\begin{split}
& \lVert \rho_+ \rVert_{L_\infty (0, T; B^{3 - 2 \slash p}_{q, p} (\Omega_+))} + \sum_{\ell = \pm} \lVert \bu_\ell \rVert_{L_\infty (0, T; B^{2 (1 - 1 \slash p)}_{q, p} (\Omega_\ell))} + \lVert h \rVert_{L_\infty (0, T; B^{3 - 1 \slash p - 1 \slash q}_{q, p} (S_R))} \\
& \qquad \le C (\CJ + \CJ_{p, q, T} (\rho_+, \bu_+, \bu_-, \pi_-, h; 0))
\end{split}
\end{equation}
holds. Hence, we see that
\begin{equation*}
E_{p, q, T} (\rho_+, \bu_+, \bu_-, \pi_-, h; 0) \le C (\CJ + \CJ_{p, q, T} (\rho_+, \bu_+, \bu_-, \pi_-, h; 0)),
\end{equation*}
where $C$ is a constant independent of $T$. Since we choose $\varepsilon$ small enough eventually, we may suppose that $0 < \varepsilon < 1$. \par
Let $T_0 > 2$ be a given number. From Theorem~\ref{th-local}, there exists a small number $\varepsilon_1 > 0$ such that the system~\eqref{eq-1.5} admits a unique solution $(\rho_+, \bu_+, \bu_-, \pi_-, h) \in \CS_{p, q} (0, 2)$ with $E_{p, q, 2} (\rho_+, \bu_+, \bu_-, \pi_-, h; 0) \le \varepsilon_1$ provided that $\CJ \le \varepsilon_1 < 1$. Assume the existence of a unique solution $(\rho_+, \bu_+, \bu_-, \pi_-, h) \in \CS_{p, q} (0, T_0)$ of \eqref{eq-1.5} satisfying 
\begin{equation}
\label{keep}
E_{p, q, T_0} (\rho_+, \bu_+, \bu_-, \pi_-, h; 0) \le \varepsilon_1, \qquad \sup_{t \in (0, T_0)} \lVert \Psi (\cdot, t) \rVert_{H^1_\infty (\dot \Omega)} \le \delta.
\end{equation}
In the following, we shall show that the solution to \eqref{eq-1.5} can be prolonged beyond $T$ keeping the estimates \eqref{keep} provided that $\varepsilon > 0$ is small enough. To this end, it suffices to show the inequality
\begin{equation}
\label{apriori}
E_{p, q, T} (\rho_+, \bu_+, \bu_-, \pi_-, h; \eta) \le M (\CJ + E_{p, q, T} (\rho_+, \bu_+, \bu_-, \pi_-, h; \eta)^2)
\end{equation}
for any $T \in (0, T_0]$ with some constant $M > 0$ independent of $\varepsilon$, $T$, and $T_0$, where $\eta$ is the same constant as in Theorem~\ref{linear-decay}. In fact, if $\CJ \le \varepsilon \ll 1$, we may deuce that $E_{p, q, T} (\rho_+, \bu_+, \bu_-, \pi_-, h; \eta) \le 2 M \varepsilon$ for any $T \in (0, T_0)$. Especially, setting $T_1 = T_0 - 1 \slash 2$, we obtain
\begin{equation*}
\lVert \rho_+ (\cdot, T_1) \rVert_{B^{3 - 2 \slash p}_{q, p} (\Omega_+)} + \sum_{\ell = \pm} \lVert \bu_\ell (\cdot, T_1) \rVert_{B^{2 (1 - 1 \slash p)}_{q, p} (\Omega_\ell)} + \lVert h (\cdot, T_1) \rVert_{B^{3 - 1 \slash p - 1 \slash q}_{q, p} (S_R)} \le 2 M \varepsilon.
\end{equation*}
Thus, choosing $\varepsilon > 0$ small enough and employing the same argument as that in proving Theorem~\ref{th-local}, we find that there exists a unique solution $(\wt \rho_+, \wt \bu_+, \wt \bu_-, \wt \pi_-, \wt h) \in \CS_{p, q} (T_1, T_1 + 1)$ of the following system:
\begin{align*}
\left\{\begin{aligned}
\pd_t \wt \rho_+ + \rho_{* +} \dv \wt \bu_+ & = f_M (\wt \rho_+, \wt \bu_+, \wt h) &\enskip &\text{ in $\Omega_+ \times (T_1, T_1 + 1)$}, \\
\rho_{* -} \dv \wt \bu_- = f_d (\wt \bu_-, \wt h) & = \rho_{* -} \dv \bff_d (\wt \bu_-, \wt h) &\enskip &\text{ in $\Omega_- \times (T_1, T_1 + 1)$}, \\
\rho_{* +} \pd_t \wt \bu_+ - \DV \bT_+ (\wt \rho_+, \wt \bu_+) & = \bff_+ (\wt \rho_+, \wt \bu_+, \wt h) &\enskip &\text{ in $\Omega_+ \times (T_1, T_1 + 1)$}, \\
\rho_{* -} \pd_t \wt \bu_- - \DV \bT_- (\wt \bu_-, \wt \pi_-) & = \bff_- (\wt \bu_-, \wt h) &\enskip &\text{ in $\Omega_- \times (T_1, T_1 + 1)$}, \\
\pd_t \wt h - \frac{1}{\rho_{* -} - \rho_{* +}} \lbrac \langle \rho_* \wt \bu, \bn \rangle \rbrac + \CM \wt \bu & = d (\wt \rho_+, \wt \bu_+, \wt \bu_-, \wt h) &\enskip &\text{ on $S_R \times (T_1, T_1 + 1)$}, \\
\bB (\wt \rho_+, \wt \bu_+, \wt \bu_-, \wt \pi_-, \wt h) & = \bG (\wt \rho_+, \wt \bu_+, \wt \bu_-, \wt h), &\enskip &\text{ on $S_R \times (T_1, T_1 + 1)$}, \\
\wt \bu_+ = 0, \qquad \langle \nabla \wt \rho_+, \bn_+ \rangle & = 0 &\enskip &\text{ on $\Gamma_+ \times (T_1, T_1 + 1)$}, \\
(\wt \rho_+, \wt \bu, \wt h) \rvert_{t = T_1} & = (\rho_+ (\cdot, T_1), \bu (\cdot, T_1), h (\cdot, T_1)) &\enskip &\text{ on $\Omega_+ \times \dot \Omega \times S_R$},	
\end{aligned}\right.	
\end{align*} 
which satisfies
\begin{equation*}
\sup_{t \in (T_1, T_1 + 1)} \lVert \Psi (\cdot, t) \rVert_{H^1_\infty (\dot \Omega)} \le \delta, \qquad  E_{p, q, (T_1, T_1 + 1)} (\wt \rho_+, \wt \bu_+, \wt \bu_-, \wt \pi_-, \wt h; 0) \le \varepsilon_1.
\end{equation*}
Here, $E_{p, q, (T_1, T_1 + 1)}$ is given by $E_{p, q, T}$ with the time interval $(T_1, T_1 + 1)$ instead of $(0, T)$. Let
\begin{alignat*}4
\ov \rho_+ & = \begin{cases}
\rho_+ & (0 < t \le T_1), \\
\wt \rho_+ & (T_1 < t < T_1 + 1),
\end{cases}
& \qquad 
\ov \bu & = \begin{cases}
\bu & (0 < t \le T_1), \\
\wt \bu & (T_1 < t < T_1 + 1),
\end{cases} \\
\ov \pi_- & = \begin{cases}
\pi_- & (0 < t \le T_1), \\
\wt \pi_- & (T_1 < t < T_1 + 1),
\end{cases}
& \qquad
\ov h & = \begin{cases}
h & (0 < t \le T_1), \\
\wt h & (T_1 < t < T_1 + 1),
\end{cases}
\end{alignat*}
and then $(\ov \rho_+, \ov \bu_+, \ov \bu_-, \ov \pi_-, \ov h)$ belongs to $\CS_{p, q} (0, T_1 + 1)$ that satisfies
\begin{equation*}
\sup_{t \in (0, T_1 + 1)} \lVert \Psi (\cdot, t) \rVert_{H^1_\infty (\dot \Omega)} \le \delta, \qquad E_{p, q, T_1 + 1} (\ov \rho_+, \ov \bu_+, \ov \bu_-, \ov \pi_-, \ov h; 0) \le \varepsilon_1 
\end{equation*}
and the system \eqref{eq-1.5} in the time interval $(0, T_1 + 1)$ instead of $(0, T)$. Since $T_1 + 1 = T_0 + 1 \slash 2$, repeating the above argument, we can prolong the solution to time interval $(0, \infty)$. This completes the proof of Theorem~\ref{th-main}. \par
Below, we show the a priori estimate \eqref{apriori}. Since we will choose $\CJ$ small enough eventually, we may assume that $0 < \CJ < 1$. We extend the right-hand members $(f_M, f_d, \bff_d, \bff_+, \bff_-, d, g, f^+_B, f^-_B, \bh, k_-)$ of \eqref{eq-1.5} to $t \in \BR$ that are denoted by $(\ov f_M, \ov f_d, \ov \bff_d, \ov \bff_+, \ov \bff_-, \ov d, \ov g, \ov f^+_B, \ov f^-_B, \ov \bh, \ov k_-)$. Here, we refer to \cite[Sec.~7.2]{W20} for the suitable extension operators. Then, employing the same argument as in \cite[Sec.~7.2]{W20}, we have
\begin{equation*}
\label{nonlinear-est}
\begin{split}
& \lVert e^{\eta t} \ov f_M (\rho_+, \bu_+, h) \rVert_{L_p (\BR; H^1_q (\Omega_+))} + \lVert e^{\eta t} \ov f_d (\bu_-, h) \rVert_{L_p (\BR; H^1_q (\Omega_-))} + \lVert e^{\eta t} \ov \bff_d (\bu_-, h) \rVert_{H^{1 \slash 2}_p (\BR; L_q (\Omega_-))} \\
& \quad + \lVert e^{\eta t} \pd_t \ov \bff_d (\bu_-, h) \rVert_{L_p (\BR; L_q (\Omega_-))} +  \lVert e^{\eta t} \ov \bff_+ (\rho_+, \bu_+, h) \rVert_{L_p (\BR; L_q (\Omega_+))} + \lVert e^{\eta t} \ov \bff_- (\bu_-, h) \rVert_{L_p (\BR; L_q (\Omega_-))} \\
& \quad + \lVert e^{\eta t} (\ov g (\rho_+, \bu_+, \bu_-, h), \nabla \ov \bh (\bu_+, \bu_-, h)) \rVert_{H^{1 \slash 2}_p (\BR; L_q (\dot \Omega))} \\
& \quad + \lVert e^{\eta t} (\nabla \ov g (\rho_+, \bu_+, \bu_-, h), \pd_t \ov \bh (\bu_+, \bu_-, h), \nabla \ov \bh (\bu_+, \bu_-, h)) \rVert_{L_p (\BR; L_q (\dot \Omega))} \\
& \quad + \sum_{\ell = \pm} \Big(\lVert e^{\eta t} \ov f^\ell_B (\rho_+, \bu_+, \bu_-, h) \rVert_{L_p (\BR; L_q (\dot \Omega))} + \lVert e^{\eta t} \ov f^\ell_B (\rho_+, \bu_+, \bu_-, h) \rVert_{H^{1 \slash 2}_p (\BR; L_q (\dot \Omega))} \\
& \qquad + \lVert e^{\eta t} \nabla \ov f^\ell_B (\rho_+, \bu_+, \bu_-, h) \rVert_{L_p (\BR; L_q (\dot \Omega))} \Big) + \lVert e^{\eta t} \nabla \ov k_- (\rho_+, h) \rVert_{H^{1 \slash 2}_p (\BR; L_q (\dot \Omega))} \\
& \quad  + \lVert e^{\eta t} (\pd_t \ov k_- (\rho_+, h), \nabla^2 \ov k_-(\rho_+, h)) \rVert_{L_p (\BR; L_q (\dot \Omega))} \\
& \le C (\CJ + E_{p, q, T} (\rho_+, \bu_+, \bu_-, \pi_-, h; \eta)^2).
\end{split}
\end{equation*}
Hence, by Theorem \ref{linear-decay}, we obtain
\begin{align*}
& \CJ_{p, q, T} (\rho_+, \bu_+, \bu_-, \pi_-, h; \eta) \\
& \le C \bigg\{\CJ + E_{p, q, T} (\rho_+, \bu_+, \bu_-, \pi_-, h; \eta)^2 + \bigg(\int_0^T (e^{\eta s} \lvert (h (\cdot, s), \varphi_0)_{S_R} \rvert)^p \ds\bigg)^{1 \slash p} \bigg\}.
\end{align*}
Noting that
\begin{equation*}
\CJ_{p, q, T} (\rho_+, \bu_+, \bu_-, \pi_-, h; 0) \le \CJ_{p, q, T} (\rho_+, \bu_+, \bu_-, \pi_-, h; \eta),
\end{equation*}
it holds
\begin{equation*}
\begin{split}
& \lVert \rho_+ \rVert_{L_\infty (0, T; B^{3 - 2 \slash p}_{q, p} (\Omega_+))} + \sum_{\ell = \pm} \lVert \bu_\ell \rVert_{L_\infty (0, T; B^{2 (1 - 1 \slash p)}_{q, p} (\Omega_\ell))} + \lVert h \rVert_{L_\infty (0, T; B^{3 - 1 \slash p - 1 \slash q}_{q, p} (S_R))} \\
& \qquad \le C (\CJ + \CJ_{p, q, T} (\rho_+, \bu_+, \bu_-, \pi_-, h; \eta))
\end{split}
\end{equation*}
instead of \eqref{402}. Hence, it follows that
\begin{align*}
& E_{p, q, T} (\rho_+, \bu_+, \bu_-, \pi_-, h; \eta) \\
& \le C \bigg\{\CJ + E_{p, q, T} (\rho_+, \bu_+, \bu_-, \pi_-, h; \eta)^2 + \bigg(\int_0^T (e^{\eta s} \lvert (h (\cdot, s), \varphi_0)_{S_R} \rvert)^p \ds\bigg)^{1 \slash p} \bigg\}.
\end{align*}
Recalling \eqref{formula-h1}, we have
\begin{equation*}
\bigg(\int_0^T (e^{\eta s} \lvert (h (\cdot, s), \varphi_0)_{S_R} \rvert)^p \ds\bigg)^{1 \slash p} \le C E_{p, q, T} (\rho_+, \bu_+, \bu_-, \pi_-, h; \eta)^2 
\end{equation*}
so that we arrive at the estimate \eqref{apriori}. This completes the proof of Theorem~\ref{th-main}.




\end{document}